\documentclass[11pt]{elsart}
\usepackage{amsmath,amsfonts,amssymb,enumerate}
\usepackage{latexsym}
\journal{Journal of Functional Analysis}
\begin{document}
\begin{frontmatter}
%\topmargin=1mm % comment out if want exact top margin
%%\advance\topmargin 19mm
%\oddsidemargin=4.6mm
%%\advance\oddsidemargin 10mm % comment out if want exact left margin
%\evensidemargin=\oddsidemargin \textwidth=160mm
%%\advance\textwidth by-60mm
%\textheight=200mm
%%\advance\textheight by-60mm
%\footskip=20mm
%%\footheight=5mm
%\headheight=10mm
%%\headsep=0mm
%%\topskip=0mm
%%\overfullrule=0pt

%\usepackage{pstricks,pstcol,pst-node,pst-plot}
%\usepackage[dvips]{graphics}

%\calclayout

\newtheorem{theorem}{Theorem}[section]
\newtheorem{lemma}[theorem]{Lemma}
\newtheorem{corollary}[theorem]{Corollary}
\newtheorem{remark}[theorem]{Remark}
\newtheorem{example}[theorem]{Example}
\newtheorem{proposition}[theorem]{Proposition}
\newtheorem{definition}[theorem]{Definition}
\renewcommand{\epsilon}{\varepsilon}
\newcommand{\vs}{\smallskip}
\newcommand{\vv}{\smallskip \smallskip}
\newcommand{\BZ}{\mathbb Z}
\def\ve{\varepsilon}
\def\exp{{\rm exp}}
\def\res{{\rm Res}}
\def\BN{{\mathbb N}}
\def\BQ{{\mathbb Q}}
\def\BR{{\mathbb R}}
\def\BC{{\mathbb C}}
\def\BS{{\mathbb S}}
\def\CA{{\mathcal A}}
\def\CB{{\mathcal B}}
\def\CC{{\mathcal C}}
\def\CS{{\mathcal S}}
\def\CF{{\mathcal F}}
\def\CP{{\mathcal P}}
\def\FA{{\frak A}}
\def\FC{{\frak C}}
\def\FS{{\frak S}}
\def\FF{{\frak F}}
\def\tC{{\tilde C}}
\def\hC{{\hat C}}
\def\bC{{\bar C}}
\def\lv{\left|}
\def\rv{\right|}
\def\lp{\left(}
\def\rp{\right)}
\def\erf{{\rm erf}}
\def\v{\smallskip}
\numberwithin{equation}{section} \numberwithin{theorem}{section}

\title {Stability of Localized Operators}
%{Stability of some operators of non-convolution type}
%{ Equivalence of $p$-stability
% of localized infinite matrices,
%synthesis operators,
%%Riesz basis
% and integral operators for different $p$}

\author[Sogang]{Chang Eon Shin},
\ead{shinc@sogang.ac.kr}
\author[UCF]{Qiyu Sun\corauthref{cor}}
\corauth[cor]{Corresponding author.}
\ead{qsun@mail.ucf.edu}

\address[Sogang]{Department of Mathematics, Sogang University, Seoul 121-742, South
Korea}
\address[UCF]{Department of Mathematics, University of Central Florida, Orlando, FL 32816, USA}
%\subjclass[2000]{94A20, 30D10}

\begin{abstract}
Let $\ell^p, 1\le p\le \infty$, be the space of all $p$-summable sequences and  $C_a$ be the
 convolution operator  associated
with a summable sequence $a$. It is known that the $\ell^p$-
stability of the convolution operator $C_a$ for different $1\le
p\le \infty$ are equivalent to each other, i.e., if $C_a$ has $\ell^p$-stability
  for some $1\le p\le \infty$ then $C_a$ has $\ell^q$-stability
  for all $1\le q\le \infty$. In  the study of spline
approximation, wavelet analysis, time-frequency analysis, and
sampling,
 there are many localized operators of non-convolution type  whose stability is one of the basic assumptions.
 In this paper, we  consider  the stability of those localized operators including
 infinite matrices in the Sj\"ostrand class,
 synthesis operators with generating functions enveloped by  shifts of a function in
 the Wiener amalgam space, and integral operators with kernels having
  certain regularity and  decay at infinity. We show that
 the $\ell^p$- stability  (or $L^p$-stability) of  those three classes of  localized  operators
  are equivalent to each other, and we  also prove that the left inverse of those localized operators are well localized.
  \end{abstract}

\end{frontmatter}

%%==================================================================
\section{Introduction}
%%==================================================================

Given a summable sequence $a=(a(j))_{j\in {\mathbb Z}}$, the
{\em convolution operator} $C_a$ associated with the  sequence $a$
is defined by
\begin{equation}\label{convolution.def}
 C_a:\ \ell^p\ni \left(b(j)\right)_{j\in {\mathbb Z}} %=:b
\longmapsto %C_ab:=
\left(\sum_{k\in {\mathbb Z}} a(j-k)
b(k)\right)_{j\in {\mathbb Z}}\in \ell^p\end{equation} where
$\ell^p$ is the set of all $p$-summable sequences with standard norm $\|\cdot\|_{\ell^p}$.
The convolution operator $C_a$ is a bounded operator on
$\ell^p$ for any $1\le p\le \infty$ and the corresponding operator norm is bounded by the
$\ell^1$ norm of the sequence $a$.

An operator $T$ on $\ell^p$ is said to have {\em $\ell^p$-stability} if there exists a positive constant $C$ such that
\begin{equation}\label{lpstability.def0}
C^{-1} \|c\|_{\ell^p}\le  \|Tc\|_{\ell^p}\le C \|c\|_{\ell^p}\quad {\rm for \ all} \ c\in \ell^p.
\end{equation}
For the convolution operator $C_a$ associated with a summable
sequence $a=(a(j))_{j\in {\mathbb Z}}$, it is known that $C_a$ has
$\ell^p$-stability for some $1\le p\le \infty$ if and only if
\begin{equation}\label{hata.chara}
\hat a(\xi)\ne 0  \quad {\rm for \ all} \ \ \xi\in {\mathbb R}
\end{equation}
where $\hat a(\xi):=\sum_{j\in {\mathbb Z}} a(j) e^{-ij\xi}$ (c.f.
\cite{akramst01,  jiamicchelli}). Therefore the $\ell^p$-stability
of the convolution operator associated with a summable sequence
are equivalent to each other for different $1\le p\le \infty$.

%For the inverse of the
%convolution operator $C_a$, it follows from the classical Wiener's
%lemma (\cite{wiener}) that if $C_a$ has bounded
%inverse in $\ell^2$ then % its inverse is a convolution operator
%%associated with a summable sequence and hence)
%it has bounded inverse in $\ell^p$ for all $1\le p\le \infty$. The converse of the above implication  is true.

\begin{theorem}\label{convolution.equivalence.tm}
Let  $a$ be a summable sequence, and $C_a$
be the convolution operator associated with the sequence $a$. If
$C_a$ has $\ell^p$-stability for some $1\le p\le \infty$, then it has $\ell^q$-stability
for all $1\le q\le \infty$.
\end{theorem}

An equivalent formulation of the above result is that the spectrum
$\sigma_p(C_a)$ of the convolution operator $C_a$ associated with
a summable sequence $a$ as an operator on $\ell^p$ is independent
of $1\le p\le \infty$,
\begin{equation}
\sigma_p(C_a)=\sigma_q(C_a)\quad {\rm for\ all} \ 1\le p, q\le
\infty
\end{equation}
see \cite{barnes90, hulanicki, pytlik} and references therein for
the discussion on spectrum of  convolution operators.

%The above equivalence of $\ell^p$-stability for the convolution operator $C_a$
% for  different $p$ can be  proved by the dual property between the sequence spaces $\ell^p$ and $\ell^{p/(p-1)}$
%  and
%  by the interpolation theorem for linear operators.
% In fact, if $C_a$ has bounded inverse on $\ell^p$, then $C_a$ has
%bounded inverse on $\ell^{p/(p-1)}$ because of
% the dual property between sequence spaces
%$\ell^p$ and $\ell^{p/(p-1)}$ and the fact that the conjugate of the convolution
%operator $C_a$ is  the  convolution operator $C_{\tilde a}$
%associated  with the sequence $\tilde a=(a(-j))_{j\in {\mathbb
%Z}}$. Therefore $C_a$ has bounded inverse on $\ell^2$ since
%$\ell^2$ is an interpolation space between $\ell^p$ and
%$\ell^{p/(p-1)}$, and the conclusion then follows from the classical Wiener's lemma.

\bigskip
Inspired by the
commutator technique developed in \cite{sjostrand} and norm
equivalence technique for a finite-dimensional space in
\cite{ald-sl}, we will give a new  proof of  Theorem \ref{convolution.equivalence.tm} in this paper
without using the characterization \eqref{hata.chara}. More importantly  we can extend the
 equivalence for $\ell^p$-stability in Theorem  \ref{convolution.equivalence.tm}
   to various  localized operators of non-convolution type, that arise in
the study of spline approximation (\cite{deboor76,demko}), wavelet
and affine frames (\cite{chuihestockler04b, jaffard}),
 Gabor frame and non-uniform sampling (\cite{balanchl04, grochenigl03, grocheniglappear,
 sunsiam}), and pseudo-differential operators (\cite{boulkhemair1,
boulkhemair2, farellstrohmer, grochenigappear, grochenigstrohmer,
herau, sjostrand, sjostrand2, toft1, toft2, toft3}).

\bigskip

Denote  by ${\mathcal A}$  the  {\em Schur class} of  infinite matrices
$A=(a(j,j'))_{j,j'\in {\mathbb Z}}$  such that
\begin{equation}\|A\|_{\mathcal A}:=\max\left(\sup_{j\in {\mathbb Z}}
\sum_{j'\in {\mathbb Z}} |a(j,j')|,\  \sup_{j'\in {\mathbb Z}}
\sum_{j\in {\mathbb Z}} |a(j,j')|\right)<\infty\end{equation}
(\cite{grocheniglappear, suntams}). An infinite matrix
$A=(a(j,j'))_{j,j'\in {\mathbb Z}}$ in the Schur class defines a
bounded operator on $\ell^p, 1\le p\le \infty$, as follows:
\begin{equation}
A:\ \ell^p\ni (c(j))_{j\in {\mathbb Z}} %=:c
\longmapsto
%Ac:=
\left(\sum_{j'\in {\mathbb Z}} a(j,j') c(j')\right)_{j\in
{\mathbb Z}}\in \ell^p.
\end{equation}
It is known that an infinite matrix $A$ is a bounded operator on
$\ell^p$ for any $1\le p\le \infty$ if and only if $A$ is in the
Schur class.
 In  Section
\ref{matrix.section}, we consider the equivalence of
$\ell^p$-stability for   infinite matrices in the Sj\"ostrand
class (see Section \ref{matrix.section} for its definition), a
subset of the Schur class ${\mathcal A}$, for different $1\le p\le
\infty$.

\begin{theorem}\label{infinitematrixpreliminary.thm}
Let $A=(a(j,j'))_{j,j'\in {\mathbb Z}}$ be an  infinite matrix
with the property that $\sum_{k\in {\mathbb Z}} \sup_{j\in
{\mathbb Z}} |a(j, j+k)|<\infty$. If $A$ has $\ell^p$-stability
for some $1\le p\le \infty$, then $A$ has $\ell^q$-stability for
all $1\le q\le\infty$.
\end{theorem}

The result in Theorem \ref{infinitematrixpreliminary.thm} for
$p=2$ follows from the Wiener's lemma in  \cite{sjostrand}. The
equivalence of  $\ell^q$-stability of an infinite matrix
$A=(a(j,j'))_{j,j'\in {\mathbb Z}}$ with $\sum_{k\in {\mathbb Z}}
\sup_{j\in {\mathbb Z}} |a(j, j+k)| (1+|k|)^s<\infty$ is
established in \cite{ald-sl} for $s>4$ and later improved in
\cite{tessera} for $s>0$.

We observe that a convolution operator $C_a$ associated with a
summable sequence $a=(a(j))_{j\in {\mathbb Z}}$ is  the operator
associated with the infinite matrix $A=(a(j-j'))_{j,j'\in {\mathbb
Z}}$ in the Sj\"ostrand class. Therefore Theorem
\ref{convolution.equivalence.tm}  follows from Theorem
\ref{infinitematrixpreliminary.thm}.

We conjecture that the equivalence of $\ell^p$-stability for
different $p\in [1, \infty]$ holds for any infinite matrix in the
Schur class ${\mathcal A}$. Some progress on the above conjecture
is made in \cite{tessera} under additional assumption that the
infinite matrix has  rows supported in balls of bounded radii.

\bigskip
For a continuous function $f$ on ${\mathbb R}$, we define the
{\em modulus of continuity } $\omega_\delta(f)$ by
\begin{equation}
\omega_\delta(f)(x)=\sup_{|y|\le \delta} |f(x+y)-f(x)|.
\end{equation}
The modulus of continuity is a delicate tool  in
mathematical analysis to measure the regularity of a function
(\cite{devore, triebel}).

For $1\le p\le \infty$, let $L^p$ be the space of all $p$-integrable functions on
${\mathbb R}^d$ with standard norm $\|\cdot\|_p$,
$${\mathcal L}^p=\left\{\phi: \ \|\phi\|_{{\mathcal L}^p}
:=\left\|\sum_{j\in{\mathbb Z}}|\phi(\cdot-j)|\right\|_{L^p([0,1])}<\infty\right\}$$
consist of functions  that are ``globally'' in  $\ell^1$ and ``locally" in $L^p$ (\cite{jiamicchelli}),
and
let
$${\mathcal W}_1=\left\{ \phi: \  \|\phi\|_{{\mathcal W}_1}:=
\sum_{k\in {\mathbb Z}} \sup_{x\in [0,1)} |\phi(k+x)|<\infty  \right\}$$
be the {\em Wiener amalgam  space}  that  consists of functions
that are ``locally''  in $L^\infty$
 and ``globally" in $\ell^1$ (\cite{akramgrochenig01}).
  We have the following inclusions for the above three classes of function spaces:
  \begin{equation}
 {\mathcal W}_1\subset {\mathcal L}^\infty \subset  {\mathcal L}^p \subset {\mathcal L}^1
 \end{equation}
 and
  \begin{equation}
 {\mathcal L}^p \subset  L^p
 \end{equation}
 where $1\le p\le \infty$ (\cite{akramgrochenig01, jiamicchelli}).

 For a family of functions
$\Phi=\{\phi_j\}_{j\in {\mathbb Z}}$ enveloped by a function $h\in
{\mathcal L}^p$, i.e.,
$$|\phi_j(x)|\le
h(x-j) \ \ {\rm for \ all} \ j\in {\mathbb Z} \ {\rm and} \ x\in
{\mathbb R},$$
 the {\em synthesis operator} $S_\Phi$  associated with $\Phi$,
\begin{equation}\label{sphi.def0}
S_\Phi:  \ell^p\ni (c(j))_{j\in {\mathbb Z}}%:c
\longmapsto %S_\Phi c:=
\sum_{j\in {\mathbb Z}} c(j) \phi_j\in L^p,
\end{equation}
is a bounded operator from $\ell^p$ to $L^p$ (\cite{jiamicchelli,
sunaicm}). For the family of functions $\Phi$ generated by  shifts
of  finitely many functions, that is,
$\Phi=\{\phi_n(\cdot-j)\}_{1\le n\le N, j\in {\mathbb Z}}$  for some
$\phi_n\in {\mathcal L}^\infty, 1\le n\le N$, it is proved in
\cite{jiamicchelli} that if the synthesis operator $S_\Phi$ in
\eqref{sphi.def0} has $\ell^p$-stability for some $1\le p\le
\infty$, i.e., there exists a positive constant $C$ such that
 \begin{equation}
 C^{-1} \|c\|_{\ell^p}\le \|S_\Phi c\|_p\le C \|c\|_{\ell^p} \quad {\rm for \ all} \ c\in \ell^p,\end{equation}
 then it has  $\ell^q$-stability  for all $1\le q\le \infty$. In Section \ref{priesz.section},
  we establish the equivalence of the stability of the
synthesis operator $S_\Phi$  for $1\le p\le
\infty$ under some regularity and decay assumption on the generating family
$\Phi$ of functions.

\begin{theorem}\label{prieszpreliminary.thm} Let
$\Phi=\{\phi_j\}_{j\in
 {\mathbb Z}} $  be the family of functions  with the property that
 \begin{equation} \label{prieszpreliminary.thm.eq1}  |\phi_j(x)|\le h(x-j)\  \ {\rm for \ all} \ x\in {\mathbb R}, j\in {\mathbb Z}, \end{equation}
and
\begin{equation}\label{prieszpreliminary.thm.eq2}
  |\omega_\delta(\phi_j)(x)|\le \omega(\delta) h(x-j)
 \ {\rm for \ all} \ x\in {\mathbb R}, j\in {\mathbb Z} \ {\rm and} \ \delta\in (0,1),\end{equation}
 where $h$ is a continuous function  in the Wiener amalgam  space ${\mathcal W}_1$ and
  $\omega(\delta)$ is a positive increasing function  on $(0,1)$ with  $\lim_{\delta\to 0} \omega(\delta)=0$,
Let $S_\Phi$ be the synthesis operator associated with the family $\Phi$ of functions. If the synthesis operator $S_\Phi$ has
  $\ell^p$-stability  for some $1\le p\le \infty$, then it
 has $\ell^q$-stability for all $1\le q\le \infty$.
\end{theorem}

For $p=2$, the conclusion  in Theorem \ref{prieszpreliminary.thm} follows from the Wiener's lemma
  in \cite{sunaicm}  under the decay assumption  (\ref{prieszpreliminary.thm.eq1}) on the
   generating family $\Phi$, but without the regularity assumption (\ref{prieszpreliminary.thm.eq2})
on the generating family $\Phi$ as in Theorem \ref{prieszpreliminary.thm}.
% We conjecture that the conclusion in Theorem \ref{prieszpreliminary.thm} holds when the function $h$ in Theorem \ref{prieszpreliminary.thm}
% is in the space ${\mathcal L}^\infty$ instead of  in the Wiener amalgam space ${\mathcal W}_1$.

\bigskip

In Section \ref{localizedintegraloperator.section}, we consider the  $L^p$-stability  of
 localized integral operators
for different $p$ (\cite{Kurbatov, sunacha}). Here a bounded
operator $S$ on $L^p$ is said to have {\em $L^p$-stability} if
there exists a positive constant $C$ such that
\begin{equation}
C^{-1} \|f\|_p\le \|Sf\|_p\le C \|f\|_p \quad {\rm for \ all} \ f\in L^p.
\end{equation}

\begin{theorem}\label{liopreliminary.thm}
Let  $I$ be the identity operator on $L^p$, and $T$ be a localized
integral operator
$$Tf(x):=\int_{{\mathbb R}} K_T(x, y) f(y)dy, \ f\in L^p$$
such that its  integral kernel $K_T$ satisfies
 \begin{equation} \label{liopreliminary.thm.eq1}
 |K_T(x,y)|\le h(x-y)\  \ {\rm for \ all} \ x, y\in {\mathbb R}, \end{equation}
and
\begin{equation}\label{liopreliminary.thm.eq2}
  |\omega_\delta(K_T)(x,y)|\le \delta^\alpha h(x-y)
 \ {\rm for \ all} \ x, y\in {\mathbb R},  {\rm and} \ \delta\in (0,1)\end{equation}
 where  $\alpha\in (0,1)$ and $h$ is a continuous function  in the Wiener amalgam  space ${\mathcal W}_1$.
If $I+T$ has $L^p$-stability for some $1\le p\le \infty$, then $I+T$ has  $L^q$-stability for all $1\le q\le \infty$.
\end{theorem}

For $p=2$, the conclusion  in Theorem \ref{liopreliminary.thm} follows from the Wiener's lemma
  in \cite{sunacha}.

\bigskip

  In this paper, we will use the following notation:
 Denote by $\chi_E$ the characteristic function on the set $E$.
 For a discrete set $\Lambda$ and $1\le p\le \infty$,  denote by $\ell^p(\Lambda)$ the set of all $p$-summable sequence
 $c:=(c(\lambda))_{\lambda\in \Lambda}$ with standard norm $\|\cdot\|_{\ell^p(\Lambda)}$, and
 by  $(\ell^p(\Lambda))^N$  the $N$ copies of $\ell^p(\Lambda)$ with norm  $\|\cdot\|_{(\ell^p(\Lambda))^N}$.
Let  ${\mathcal B}(\ell^p(\Lambda)), 1\le p\le \infty$,  be
  the Banach algebra containing
 all bounded operators on $\ell^p(\Lambda)$ embedded with  standard operator norm.
 For a measurable set $E$, let $L^p(E)$ be the space of all $p$-integrable functions on $E$ with norm $\|\cdot\|_{L^p(E)}$, while if $E={\mathbb R}^d$,  we use $L^p$ instead of $L^p({\mathbb R}^d)$ and
$\|\cdot\|_p$ instead of $\|\cdot\|_{L^p({\mathbb R}^d)}$ for brevity.
Let ${\mathcal B}(L^p)$ be
  the Banach algebra containing
 all bounded operators on $L^p({\mathbb R}^d)$ embedded with  standard operator norm.
For $x=(x_1, \ldots,
x_d)\in
{\mathbb R}^d$, define
$\|x\|_\infty=\max(|x_1|, \ldots, |x_d|)$.
Denote by $I$ the identity operator on $L^p, 1\le p\le \infty$, or  the identity matrix of appropriate size.

  In this paper, the uppercase letter $C$ denotes an absolute
constant that may be different at different occurrences, except
stated explicitly.

%%==================================================================
\section{$\ell^p$-stability for localized  infinite matrices}
\label{matrix.section}
%%==================================================================

In this section,  we  consider the $\ell^p$-stability for infinite
matrices of the form $\big(a(\lambda,\lambda') \big)_{\lambda\in
\Lambda,\lambda'\in \Lambda'}$ having certain off-diagonal decay.
That kind of extreme noncommutative matrices  arises in the study
of spline approximation (\cite{deboor76,demko}), wavelet and
affine wavelets (\cite{chuihestockler04b, jaffard}), Gabor frame
(\cite{balanchl04, grochenigl03, grocheniglappear}), non-uniform
sampling (\cite{sunsiam}), and pseudo-differential operators
(\cite{boulkhemair1, boulkhemair2, farellstrohmer,
grochenigappear, grochenigstrohmer, herau, sjostrand, sjostrand2,
toft1, toft2, toft3}), and the $\ell^p$-stability for those
matrices is one of few basic assumptions in these studies. The
main results of this section are Theorem \ref{matrixstability.tm}
(a slight generalization of Theorem
\ref{infinitematrixpreliminary.thm}), Corollary
\ref{slantmatrix.cor} (an equality for spectrum of slanted
matrices on $\ell^p$ for different $p$) and Corollary
\ref{wienerinfinitematrix.cr} (a Wiener's lemma for the
Sj\"ostrand class of infinite matrices).

To state our result on stability for localized  infinite matrices,
we recall three concepts. We say that a discrete subset $\Lambda$
of ${\mathbb R}^d$  is  {\em relatively-separated } if
\begin{equation}\label{relativelyseparated.def}
R(\Lambda)=\sup_{x\in {\mathbb R}^d} \sum_{\lambda\in \Lambda}
\chi_{\lambda+[0,1)^d}(x)<\infty.
\end{equation}
For relatively-separated subsets  $\Lambda, \Lambda'$  of
${{\mathbb R}^d}$, we let
${\mathcal C}(\Lambda, \Lambda')$, or ${\mathcal C}$ for short,
 be the {\em Sj\"ostrand  class} of infinite matrices $A=(a(\lambda, \lambda'))_{\lambda\in \Lambda,
\lambda'\in \Lambda'}$ such that
\begin{equation}\label{sjostrand.def}
\|A\|_{{\mathcal C}}:=\sum_{k\in {\mathbb Z}^d}
 \sup_{\lambda\in \Lambda, \lambda'\in \Lambda'} |a(\lambda, \lambda')| \chi_{k+[0,1)^d}(\lambda-\lambda').
\end{equation}
%where
%$$a(k)=\left\{\begin{array}{ll} \sup_{\lambda-\lambda'\in k+[0,1)^d}|a(\lambda, \lambda')| & {\rm  if}\ (\Lambda-\Lambda')\cap (k+[0,1)^d)\ne \emptyset\\
%0 &  {\rm if}  \ (\Lambda-\Lambda')\cap (k+[0,1)^d)=\emptyset\end{array}\right.$$
(\cite{sjostrand, suntams}). As usual, an infinite matrix
$A=(a(\lambda, \lambda'))_{\lambda\in \Lambda, \lambda'\in
\Lambda'}$ in the Sj\"ostrand class ${\mathcal C}(\Lambda,
\Lambda')$ defines a bounded operator from $\ell^p(\Lambda')$ to
$\ell^p(\Lambda)$,
\begin{equation}
A:\ \ell^p(\Lambda')\ni (c(\lambda'))_{\lambda'\in
\Lambda'}:=c\longmapsto Ac:= \left( \sum_{\lambda'\in \Lambda'}
a(\lambda, \lambda') c(\lambda')\right)_{\lambda\in \Lambda}\in
\ell^p(\Lambda),
\end{equation}
where $1\le p\le \infty$ (\cite{suntams}). For $1\le p\le \infty$,
we say that an infinite matrix $A:=(a(\lambda,
\lambda'))_{\lambda\in \Lambda, \lambda'\in \Lambda'}$ has {\em
$\ell^p$-stability} if there exists a positive constant $C$ such
that
\begin{equation}
C^{-1} \|c\|_{\ell^p(\Lambda')} \le \|Ac\|_{\ell^p(\Lambda)} \le C
\|c\|_{\ell^p(\Lambda')} \quad {\rm for \ all} \ c\in
\ell^p(\Lambda').
\end{equation}

\begin{theorem}\label{matrixstability.tm}
Let $\Lambda, \Lambda'$ be two relatively-separated subsets of
${\mathbb R}^d$, and $A=(a(\lambda, \lambda'))_{\lambda\in
\Lambda,\lambda'\in \Lambda'}$ be an  infinite matrix in the
Sj\"ostrand class ${\mathcal C}(\Lambda, \Lambda')$. If $A$ has
$\ell^p$-stability for some $1\le p\le \infty$, then $A$ has
$\ell^q$-stability for all $1\le q\le \infty$.
\end{theorem}

For $\Lambda'=\Lambda$, the above theorem can be reformulated as
follows:

\begin{corollary} Let $\Lambda$ be a relatively-separated subset
of ${\mathbb R}^d$. Then the spectrum $\sigma_p(A)$ of an infinite
matrix $A$ in the Sj\"ostrand class ${\mathcal C}(\Lambda,
\Lambda)$  as an operator on $\ell^p(\Lambda)$ is independent of
$1\le p\le \infty$, i.e.,
\begin{equation}
\sigma_p(A)=\sigma_q(A) \quad \ 1\le p, q\le \infty.
\end{equation}
%holds for any $A\in {\mathcal C}(\Lambda, \Lambda)$ and $1\le p,
%q\le \infty$.
\end{corollary}

\bigskip

  A function $w$ on ${\mathbb R}^d$ is said to be a
 {\em weight} if $w(x)\ge 1$ for all $x\in {\mathbb R}^d$ and
 $\sup_{|y|\le 1}\sup_{x\in {\mathbb R}^d}\frac{w(x  +y)}{w(x)}<\infty$.
  For a weight $w$ on ${\mathbb R}^d$ and a positive number $\alpha$,  denote by $\Sigma_\alpha^w$
  the  family of all {\em $\alpha$-slant infinite matrices}
 $A=(a(j,j'))_{j,j'\in {\mathbb Z}^d}$ with
 $$  \|A\|_{\Sigma_\alpha^w} =\sum_{k\in{\mathbb Z}^d} w(k) \sup_{j, j'\in{\mathbb Z}^d} |a(j, j')|
\chi_{k+[0,1)^d}(j'-\alpha j)<\infty.
$$
The slanted matrices appear  in wavelet theory, signal processing
and sampling  theory \cite{ald-sl, berg, dm, kovacevic}, and also
occur in the  K-theory of operator algebras \cite{yu}.
 Note that
\begin{equation}
\Sigma_\alpha^w\subset {\mathcal C}(\alpha{\mathbb Z}^d, {\mathbb Z}^d)
\end{equation}
for any weight $w$. Then we obtain  the following result from
Theorem \ref{matrixstability.tm}.

\begin{corollary}\label{slantmatrix.cor}
 Let $\alpha>0$ and
 $w$ be a weight.  If $A\in \Sigma_\alpha^w$
has $\ell^p$-stability for some $1\le p\le \infty$, then $A$ has
$\ell^q$-stability for all $1\le q\le \infty$.
\end{corollary}

The result in the above corollary is established in \cite{ald-sl}
for the weight function $w(x)=(1+|x|)^s$ with $s>(d+1)^2$, and in
\cite{tessera} for  the weight function $w(x)=(1+|x|)^s$ with
$s>0$.

\bigskip

Given a Banach algebra ${\mathcal B}$, we say that a subalgebra
${\mathcal A}$ of ${\mathcal B}$ is {\em inverse-closed}
if  the inverse $T^{-1}$ of the operator $T\in {\mathcal A}$
belongs to ${\mathcal B}$ implies that it
 belongs to ${\mathcal A}$ (\cite{connes, gelfand, naimark, rieffel,  takesaki}).
 The inverse-closed subalgebra was first studied  for periodic
functions with  absolutely convergent Fourier series, which states
that  {\em if a periodic function $f$ does not vanish on the real
line  and has absolutely convergent Fourier series, i.e.,
 $f(x)=\sum_{j\in {\mathbb Z}} a(j) e^{-ijx}$
 and $\sum_{j\in {\mathbb Z}}
|a(j)|<\infty$,
 then    $f^{-1}$ has  absolutely convergent Fourier series too}
 (\cite{wiener}).
An equivalent formulation of the above  Wiener's lemma involving
matrix algebras is that the commutative Banach algebra
\begin{equation}\label{eq1.8}
 \tilde {\mathcal W}:=\left\{\big(a(j-j') \big)_{j,j'\in {\mathbb Z}}, \
\sum_{j\in  {\mathbb Z}} |a(j)|<\infty\right\}
\end{equation}
is an inverse-closed  Banach subalgebra of ${\mathcal
B}^2(\ell^2({\mathbb Z}))$ (\cite{wiener}).
 The classical Wiener's lemma and its various generalizations
(see, e.g., \cite{balanchl04, balank09, barnes, baskakov, connes,
demko, grochenigl03, grocheniglappear, domar56, jaffard, rieffel,
sjostrand}) are important and
 have numerous applications
  in numerical analysis, wavelet theory, frame theory, and sampling theory.
 For example, the classical Wiener's lemma
 and its weighted variation (\cite{domar56}) were
used to establish the decay property at infinity for  dual
generators of a shift-invariant space (\cite{akramgrochenig01,
jiamicchelli}); the Wiener's lemma for matrices associated with
twisted convolution was used in the study the decay properties of
the dual Gabor frame for $L^2$ (\cite{balanchl04, grochenigl03,
grocheniglappear}); the Jaffard's result (\cite{jaffard}) for
infinite matrices with polynomial decay was used in numerical
analysis (\cite{christensenstrohmer, strohmer00, strohmer02}),
wavelet analysis (\cite{jaffard}), time-frequency analysis
(\cite{grochenigbook, grochenig03, grochenig04}) and sampling
(\cite{atreasbenedetto, corderogro04, grochenig04, sunaicm});  and
the Sj\"ostrand's result (\cite{sjostrand}) for infinite matrices
was used in the study of pseudo-differential operators, Gabor
frames and sampling (\cite{balanchl04, grochenigappear, sjostrand,
sunsiam}). Therefore there are lots of papers devoted to the
Wiener's lemma for infinite matrices with various off-diagonal
decay conditions (see \cite{balan, balanchl04, balank09, baskakov,
bochner, demko,
 grocheniglappear, domar56, jaffard, sjostrand, sun1,
suntams} and also \cite{grochenigrsappear} for a short historical
review). The Wiener's lemma for
 the Sj\"ostrand class ${\mathcal C}(\Lambda, \Lambda)$ of
 infinite matrices $(a(\lambda, \lambda'))_{\lambda, \lambda'\in \Lambda}$
 says that  ${\mathcal C}(\Lambda, \Lambda)$ is an inverse-closed subalgebra of
 ${\mathcal B}(\ell^2(\Lambda))$ where $\Lambda$ is a relatively-separated subset  of ${\mathbb R}^d$ (\cite{sjostrand}).
  This together with the equivalence of $\ell^p$-stability  for different $p$ in Theorem \ref{matrixstability.tm}
  proves that ${\mathcal C}(\Lambda, \Lambda)$ is an inverse-closed subalgebra of
 ${\mathcal B}(\ell^p(\Lambda))$ for any $1\le p\le \infty$.

 \begin{corollary}\label{wienerinfinitematrix.cr}
 Let $1\le p\le \infty$ and $\Lambda$ be a relatively-separated subset of
${\mathbb R}^d$. Then the
Sj\"ostrand class ${\mathcal C}(\Lambda, \Lambda)$ is an inverse-closed subalgebra of
${\mathcal B}(\ell^p(\Lambda))$, i.e., if $A\in  {\mathcal C}(\Lambda, \Lambda)$
has bounded inverse on ${\mathcal B}(\ell^p(\Lambda))$, then $A^{-1}\in {\mathcal C}(\Lambda, \Lambda)$.
 \end{corollary}

Before we start the proof of Theorem \ref{matrixstability.tm}, let
us  consider necessary conditions on the relatively-separated
subsets $\Lambda$ and $\Lambda'$ such that there exists a matrix
$A=(a(\lambda, \lambda'))_{\lambda\in \Lambda, \lambda'\in
\Lambda'}$ in the Sj\"ostrand class ${\mathcal C}(\Lambda,
\Lambda')$ which has $\ell^p$-stability for some $1\le p\le
\infty$. Similar conclusion is obtained in \cite{sunsiam} for
sampling signals with finite rate of innovation, and in
\cite{pfander} for slanted matrices.

\begin{proposition} Let $\Lambda, \Lambda'$ be
relatively-separated subsets of ${\mathbb R}^d$. If there exists
 a matrix
$A=(a(\lambda, \lambda'))_{\lambda\in \Lambda, \lambda'\in
\Lambda'}$ in the Sj\"ostrand class ${\mathcal C}(\Lambda,
\Lambda')$ which has $\ell^p$-stability for some $1\le p\le
\infty$, then  there exists a positive number $R_0$ such that for
any bounded set $K$ the cardinality of the set $\Lambda\cap B(K,
R_0)$ is larger than or equal to the cardinality of the set
$\Lambda'\cap K$, where $B(K,R)$ is the set of all points in
${\mathbb R}^d$ with distance to $K$ less than $R$.
\end{proposition}

 \begin{pf}
 We  show
the above result on relatively-separated sets $\Lambda$ and
$\Lambda'$ by similar argument to the proof of the necessary
condition on the sampling set of a stable sampling and
reconstruction process in \cite{sunsiam}. Let $K$ be a compact
subset of ${\mathbb R}^d$ and $\ell^p( \Lambda'\cap K)$ be the space of
all sequences in $\ell^p(\Lambda')$ supported on $\Lambda'\cap K$.
For a sequence $c\in \ell^p(\Lambda'\cap K)$, it follows from the
property of the matrix  $A$ in the S\"ostrand class that the
$\ell^p$ norm of the sequence $Ac$ outside of $B(K, R)\cap
\Lambda$ is less than $\epsilon(R)\|c\|_{\ell^p(\Lambda')}$, where
$\epsilon(R)$ (independent of  the compact set $K$)  tends to zero
as $R$ tends to infinity. Thus there exists a positive constant
$R_0$ such that the $\ell^p$ norm of the sequence  $Ac$ inside
$B(K, R_0)\cap \Lambda$ is equivalent to the $\ell^p$ norm of the
sequence $c$. This implies that the submatrix obtained by
selecting the columns in $\Lambda\cap B(K, R_0)$ and rows in
$\Lambda'\cap K$  of the matrix $A$ has full rank, which proves
the desired conclusion on the relatively-separated subsets
$\Lambda$ and $\Lambda'$. $\qquad\Box$ \end{pf}

\bigskip

The proof of Theorem \ref{matrixstability.tm} is inspired by the
commutator technique developed in \cite{sjostrand} and norm
equivalence technique for a finite-dimensional space in
\cite{ald-sl}. To prove Theorem \ref{matrixstability.tm}, we need
several  lemmas. First we recall a known result about the
boundedness of infinite matrices in the Sj\"ostrand class.

\begin{lemma}\label{lpboundedness.lm} {\rm (\cite{sunaicm})}\
Let $1\le p\le \infty$,  $\Lambda$ and $\Lambda'$ be two
relatively-separated subsets of ${\mathbb R}^d$, and
$A=(a(\lambda, \lambda'))_{\lambda\in \Lambda,\lambda'\in
\Lambda'}$ be an infinite matrix in the Sj\"ostrand class
${\mathcal C}(\Lambda, \Lambda')$. Then the infinite matrix $A$ is
a bounded operator from $\ell^p(\Lambda')$ to $\ell^p(\Lambda)$.
Moreover there exists an absolute constant $C$ (that depends on $d$ and $p$
only) such that
\begin{equation}\label{eq-3.1}
\|A c\|_{\ell^p(\Lambda)}\le C  R(\Lambda)^{1/p}
R(\Lambda')^{1-1/p} \|A\|_{\mathcal C} \|c\|_{\ell^p(\Lambda')}
\quad
{\rm for \ all} \ c%=(c(\lambda'))_{\lambda'\in\Lambda'}
\in
\ell^p(\Lambda').
\end{equation}
\end{lemma}

\bigskip

Define the cut-off function
\begin{equation}\label{psi.def}\psi(x)=\min(\max(2-\|x\|_\infty, 0),
1)= \left\{\begin{array}{ll}
1 \  & {\rm if} \ \|x\|_\infty\le 1,\\
2-\|x\|_\infty & {\rm if} \ 1<\|x\|_\infty <2,\\
0 & {\rm if} \ \|x\|_\infty\ge 2.\end{array}\right.
 \end{equation} Then
\begin{equation}\label{cutoff.eq}
\left\{\begin{array}{l} 0\le \psi(x)\le 1 \ {\rm for \ all} \ x\in
{\mathbb R}^d,  \ {\rm and}\\
|\psi(x)-\psi(y)|\le \|x-y\|_\infty \ {\rm for \ all} \ x, y\in
{\mathbb R}^d.\end{array}\right.
\end{equation}

%for the real line, and
%$\psi(x_1, \ldots, x_d)= \prod_{i=1}^d \psi_0(x_i)$ for $(x_1, \ldots, x_d)\in {\mathbb R}^d$.
For $n\in {\mathbb Z}^d$ and $N\in {\mathbb N}$, define the
multiplication operator $\Psi_n^N :\ \ell^p(\Lambda)\to
\ell^p(\Lambda)$ by \begin{equation}
\label{multiplicationoperator.def} \Psi_n^N
c=\Big(\psi\big(\frac{\lambda-n}{N}\big)
c(\lambda)\Big)_{\lambda\in \Lambda} \ \  {\rm for} \
c=(c(\lambda))_{\lambda\in \Lambda}\in \ell^p(\Lambda)
\end{equation} where $\Lambda$ is a relatively-separated subset of
${\mathbb R}^d$. The multiplication operator $\psi_n^N$ can also
be thought as a diagonal matrix ${\rm
diag}(\psi((\lambda-n)/N))_{\lambda\in \Lambda}$.

For an infinite matrix $A=(a(\lambda, \lambda'))_{\lambda\in
\Lambda, \lambda'\in \Lambda'}$ and any $s\ge 0$, define the
truncation matrix
\begin{equation}\label{truncationmatrix.def}
A_s=(a_s(\lambda, \lambda'))_{\lambda\in \Lambda, \lambda'\in
\Lambda'}
\end{equation}
where $ a_s(\lambda, \lambda')= a(\lambda, \lambda')$  if
$\|\lambda-\lambda'\|_\infty<s$ and $a_s(\lambda, \lambda')=0$
otherwise. For the truncation matrices $A_s, s\ge 0$, of an
infinite matrix $A$ in the Sj\"ostrand class ${\mathcal
C}(\Lambda, \Lambda')$, we have
\begin{equation}
\label{truncation.eq} \|A-A_s\|_{\mathcal C} \ {\rm is \ a\
decreasing \ function \ with } \ \lim_{s\to +\infty}
\|A-A_s\|_{\mathcal C}=0.
\end{equation}

%Clearly we have
%the multiplication operators $\Psi_n^N$ are bounded operator on $\ell^p(\Lambda)$ with
%operator norm one,
%\begin{equation}
%\|\Psi_n^Nc\|_{\ell^p(\Lambda)}\le \|c\|_{\ell^p(\Lambda)}.
%\end{equation}

\begin{lemma}\label{multiplicationoperation.lm}
Let $1\le p, q\le \infty$,  $1\le N \in {\mathbb N}$, $\Lambda$
and $\Lambda'$ be two relatively-separated subsets of ${\mathbb
R}^d$, and $A=(a(\lambda, \lambda'))_{\lambda\in
\Lambda,\lambda'\in \Lambda'}$ be an infinite matrix in the
Sj\"ostrand class ${\mathcal C}(\Lambda, \Lambda')$. Then there exists an absolute constant $C$
(that depends on $d, p, q$ only) such that
\begin{eqnarray}\label{eq-3.6}
& & \left\| \left(\|(A_{N} \Psi_n^N - \Psi^{N}_{n} A_{
N})c\|_{\ell^p(\Lambda)}\right)_{n\in N{\mathbb Z}^d}\right\|_{\ell^q(N{\mathbb Z}^d)} \nonumber\\
 &\le & C R(\Lambda)^{1/p} R(\Lambda')^{1-1/p} \min_{0\le s\le
N}\Big( \|A-A_{{s}}\|_{\mathcal C}+\frac{s}{N} \|A\|_{\mathcal
  C}\Big)\nonumber\\
  & &\qquad \times
  \left\| \left(
  \|\Psi_n^N c\|_{\ell^p(\Lambda')}\right)_{n\in N{\mathbb Z}^d}\right\|_{\ell^q(N{\mathbb Z}^d)}
  \quad {\rm for \ all} \ c%=(c(\lambda'))_{\lambda'\in \Lambda'}
\in \ell^q(\Lambda').
\end{eqnarray}
%holds for any sequence $c%=(c(\lambda'))_{\lambda'\in \Lambda'}\in \ell^q(\Lambda')$.
\end{lemma}
\begin{pf} Observing that
$$
A_{N} \Psi_n^N - \Psi^{N}_{n} A_{N}= (A_{N} \Psi_n^N -\Psi^{N}_{n}
A_{N})\Psi_n^{6N},
$$
we obtain from  Lemma \ref{lpboundedness.lm} that
\begin{align}
\|(A_{N} \Psi_n^N-\Psi^{ N}_{ n} A_{N})c\|_{\ell^p(\Lambda)}
%  & =\|(A_{MN} \Psi_n^N-\Psi^{ N}_{n} A_{MN})P^{d(2N+M)}_n c\|_{\ell^p(\Lambda)}
%\nonumber
%\\
  & \le C  R(\Lambda)^{1/p} R(\Lambda')^{1-1/p}\nonumber\\
  & \quad \times \|A_{N} \Psi_n^N-\Psi^{ N}_{n} A_{N}\|_{\mathcal C}
    \|\Psi_n^{6N} c\|_{\ell^p(\Lambda')}
\label{eq-3.7}
\end{align}
for any $c\in \ell^p(\Lambda')$. We note  from \eqref{psi.def},
\eqref{cutoff.eq} and \eqref{truncation.eq} that
\begin{eqnarray}
\| A_{N} \Psi_n^N-\Psi^{ N}_{n} A_{N}\|_{\mathcal C}%\nonumber \\
& \le & \big\|\big(a_{\small N}(\lambda,\lambda')
    (\psi^N_n(\lambda')-\psi^{N}_{n}
     (\lambda)\big)_{\lambda\in \Lambda,\lambda'\in \Lambda'}\big\|_{\mathcal
     C}
     \nonumber \\
  &\le & \min_{0\le s\le N}\Big( \big\|\big(a_{\small N}(\lambda,\lambda')-a_{s}(\lambda,\lambda')
  \big)_{\lambda\in \Lambda,\lambda'\in \Lambda'}\big\|_{\mathcal C}
\nonumber \\
  & &  + \frac{s}{ N} \big\|\big(a_{s}(\lambda,\lambda')
  \big)_{\lambda\in \Lambda,\lambda'\in \Lambda'}\big\|_{\mathcal
  C}\Big)
\nonumber  \\
 &\le &  \min_{0\le s\le N}\Big( \|A-A_{{s}}\|_{\mathcal C}+\frac{s}{N} \|A\|_{\mathcal C}\Big).
 \label{eqnew-3.7}
\end{eqnarray}
Then we combine (\ref{eq-3.7}) and (\ref{eqnew-3.7}) to yield
\begin{eqnarray} \label{eq-3.8}
& & \|(A_{N} \Psi_n^N-\Psi^{ N}_{ n} A_{N})c\|_{\ell^p(\Lambda)}
\le
%  & =\|(A_{MN} \Psi_n^N-\Psi^{ N}_{n} A_{MN})P^{d(2N+M)}_n c\|_{\ell^p(\Lambda)}
%\nonumber
%\\
C R(\Lambda)^{1/p} R(\Lambda')^{1-1/p}  \nonumber\\
  &&\qquad\qquad\qquad \times \min_{0\le s\le N}\Big( \|A-A_{{s}}\|_{\mathcal C}+\frac{s}{N} \|A\|_{\mathcal
  C}\Big)
        \|\Psi_n^{6N} c\|_{\ell^p(\Lambda')}
\end{eqnarray}
for any $c\in \ell^p(\Lambda')$. Thus for $1\le q\le \infty$,  we
get from \eqref{psi.def}, \eqref{cutoff.eq} and \eqref{eq-3.8}
that
\begin{eqnarray*}
&  & \left\|\left( \big\|(A_{N}\Psi_n^N -
 \Psi_{n}^{N} A_{N})c\big\|_{\ell^p(\Lambda)}\right)_{n\in N{\mathbb Z}^d}\right\|_{\ell^q(N{\mathbb Z}^d)}
\\
  & \le & C R(\Lambda)^{1/p} R(\Lambda')^{1-1/p}
  \min_{0\le s\le N}\Big( \|A-A_{{s}}\|_{\mathcal C}+\frac{s}{N} \|A\|_{\mathcal
  C}\Big)\nonumber\\
  & &  \times
    \left\|\left( \|\Psi_n^{6N}c\|_{\ell^p(\Lambda')}\right)_{n\in N{\mathbb Z}^d}\right\|_{\ell^q(N{\mathbb Z}^d)}
    \\
  &\le &  C  R(\Lambda)^{1/p} R(\Lambda')^{1-1/p}
  \min_{0\le s\le N}\Big( \|A-A_{{s}}\|_{\mathcal C}+\frac{s}{N} \|A\|_{\mathcal C}\Big)
    \nonumber \\
&& \times    \sum_{j\in {\mathbb Z}^d \ {\rm with}\ \|j\|_\infty\le 6}
\left\|\left(\|
     \Psi_{n+2jN}^N c\|_{\ell^p(\Lambda')}\right)_{n\in N{\mathbb Z}^d}\right\|_{\ell^q(N{\mathbb Z}^d)}
\\
  & \le &  C  R(\Lambda)^{1/p} R(\Lambda')^{1-1/p}
  \min_{0\le s\le N}\Big( \|A-A_{{s}}\|_{\mathcal C}+\frac{s}{N} \|A\|_{\mathcal C}\Big)
   \nonumber\\
   & &  \times \left\|\left( \|\Psi_n^N c\|_{\ell^p(\Lambda')}\right)_{n\in N{\mathbb Z}^d}\right\|_{\ell^q(N{\mathbb Z}^d)}.
\end{eqnarray*}
This  proves  the estimate \eqref{eq-3.6}.
$\qquad \Box$
\end{pf}

\begin{lemma}\label{le-3.4}
Let $1\le  N\in {\mathbb N}, 1\le p, q\le \infty$, $\Lambda$ and
$\Lambda'$ be two relatively-separated subsets of ${\mathbb R}^d$,
and $A=(a(\lambda, \lambda'))_{\lambda\in \Lambda,\lambda'\in
\Lambda'}$ be an infinite matrix in the Sj\"ostrand class
${\mathcal C}(\Lambda, \Lambda')$. Then there exists a positive
constant $C$ (that depends only on $d, p, q$) such that
\begin{eqnarray}\label{le-3.4.eq1}
   \left\|\left(  \|\Psi_{ n}^{ N}
    A c\|_{\ell^p(\Lambda)}\right)_{n\in N{\mathbb Z}^d}\right\|_{\ell^q(N{\mathbb Z}^d)}
  & \le &  C R(\Lambda)^{1/p} R(\Lambda')^{1-1/p} \|A\|_{\mathcal
  C}\nonumber\\
  & &  \times
    \left\|\left( \|\Psi_n^N c\|_{\ell^p(\Lambda')}\right)_{n\in N{\mathbb Z}^d}\right\|_{\ell^q(N{\mathbb Z}^d)}\end{eqnarray}
holds for any sequence $c %=(c(\lambda'))_{\lambda'\in \Lambda'}
\in
\ell^q(\Lambda')$.
\end{lemma}

\begin{pf} By \eqref{psi.def} and \eqref{cutoff.eq}, we have that
\begin{eqnarray}\label{cutoff.partition}
 4^d\ge \sum_{k\in {\mathbb Z}^d} \chi_{[-2, 2)^d}(x-k) & \ge &
 \sum_{k\in {\mathbb Z}^d} (\psi(x-k))^2\nonumber\\
 & \ge &  \sum_{k\in {\mathbb Z}^d} \chi_{[-1, 1)^d}(x-k)= 2^d \quad  {\rm  for\ all}\  x\in {\mathbb R}^d.
  \end{eqnarray}
 Combining \eqref{cutoff.partition} and Lemma
\ref{lpboundedness.lm}, we obtain that
\begin{eqnarray} \label{le-3.4.pf.eq1}
\|\Psi_{ n}^{ N}
    A c\|_{\ell^p(\Lambda)} & \le & C  R(\Lambda)^{1/p} R(\Lambda')^{1-1/p}
  \sum_{n'\in N{\mathbb Z}^d}
   \|\Psi_{ n}^{ N}
    A \Psi_{n+n'}^N \|_{\mathcal C} \|\Psi_{n+n'}^N
    c\|_{\ell^p(\Lambda')}\nonumber\\
  & \le & C  R(\Lambda)^{1/p} R(\Lambda')^{1-1/p}\nonumber\\
  & &  \times
  \sum_{n'\in N{\mathbb Z}^d} \Big( \sum_{k\in {\mathbb Z}^d \ {\rm with } \ \|k-n'\|_\infty\le 4N} a(k)\Big)
  \|\Psi_{n+n'}^N c\|_{\ell^p(\Lambda')}
\end{eqnarray}
holds for $n\in N{\mathbb Z}^d$ and $c\in \ell^q(\Lambda')$, where
$$a(k)=\sup_{\lambda\in \Lambda, \lambda'\in \Lambda'} |a(\lambda,
\lambda')|\chi_{k+[0,1)^d}(\lambda-\lambda').$$ From
\eqref{le-3.4.pf.eq1}
  we get  that
\begin{eqnarray*}
  & &  \left\|\left(  \|\Psi_{ n}^{ N}
    A c\|_{\ell^p(\Lambda)}\right)_{n\in
    N{\mathbb Z}^d}\right\|_{\ell^q(N{\mathbb Z}^d)}\nonumber\\
%& \le  &
% C  R(\Lambda)^{1/p} R(\Lambda')^{1-1/p}\sum_{n\in N{\mathbb Z}^d} \Big(\sum_{n'\in N{\mathbb Z}^d}
%   \|\Psi_{ n}^{ N}
%    A \Psi_{n+n'}^N \|_{\mathcal C} \|\Psi_{n+n'}^N c\|_{\ell^p(\Lambda)}\Big)^q
%\Big)^{1\over q}\\
%& \le  & C R(\Lambda)^{1/p} R(\Lambda')^{1-1/p}
%  \sum_{n\in N{\mathbb Z}^d} \Big(\sum_{n'\in N{\mathbb Z}^d}
% \big(\sum_{k\in {\mathbb Z}^d \ {\rm with } \ |k-Nn'|\le 4N} a(k)\big)
% \|\Psi_{n+n'}^N c\|_{\ell^p(\Lambda)}\Big)^q
%\Big)^{1\over q}
%\\
& \le  & C R(\Lambda)^{1/p} R(\Lambda')^{1-1/p} \Big(\sum_{n'\in
N{\mathbb Z}^d}\ \ \sum_{k\in {\mathbb Z}^d \ {\rm with } \
\|k-n'\|_\infty\le 4N} a(k)\Big)
 \\
 & &  \times
   \left\|\left(
 \|\Psi_{n}^N c\|_{\ell^p(\Lambda')}\right)_{n\in N{\mathbb Z}^d}\right\|_{\ell^q(N{\mathbb Z}^d)}\\
& \le  &   C R(\Lambda)^{1/p} R(\Lambda')^{1-1/p} \|A\|_{\mathcal
C}    \left\|\left(
 \|\Psi_{n}^N c\|_{\ell^p(\Lambda')}\right)_{n\in N{\mathbb Z}^d}\right\|_{\ell^q(N{\mathbb Z}^d)}
\end{eqnarray*}
for $1\le q\le \infty$. %,  and
%\begin{eqnarray*}
%  \sup_{n\in N{\mathbb Z}^d} \|\Psi_{ n}^{ N}
%    A c\|_{\ell^p(\Lambda)}
%%& \le  &
%% C  R(\Lambda)^{1/p} R(\Lambda')^{1-1/p}\sum_{n\in N{\mathbb Z}^d} \Big(\sum_{n'\in N{\mathbb Z}^d}
%%   \|\Psi_{ n}^{ N}
%%    A \Psi_{n+n'}^N \|_{\mathcal C} \|\Psi_{n+n'}^N c\|_{\ell^p(\Lambda)}\Big)^q
%%\Big)^{1\over q}\\
%%& \le  & C R(\Lambda)^{1/p} R(\Lambda')^{1-1/p}
%%  \sum_{n\in N{\mathbb Z}^d} \Big(\sum_{n'\in N{\mathbb Z}^d}
%% \big(\sum_{k\in {\mathbb Z}^d \ {\rm with } \ |k-Nn'|\le 4N} a(k)\big)
%% \|\Psi_{n+n'}^N c\|_{\ell^p(\Lambda)}\Big)^q
%%\Big)^{1\over q}
%%\\
%& \le  & C R(\Lambda)^{1/p} R(\Lambda')^{1-1/p} \sum_{n'\in
%N{\mathbb Z}^d} \Big(\sum_{k\in {\mathbb Z}^d \ {\rm with } \ |k-Nn'|\le 4N}
%a(k)\Big)\\
%& & \qquad \times
% \sup_{n\in N{\mathbb Z}^d}
% \|\Psi_{n+n'}^N c\|_{\ell^p(\Lambda')}\\
%& \le  &   C R(\Lambda)^{1/p} R(\Lambda')^{1-1/p} \|A\|_{\mathcal
%C}  \sup_{n\in N{\mathbb Z}^d}
% \|\Psi_{n}^N c\|_{\ell^p(\Lambda')}.
%\end{eqnarray*}
Then the estimate \eqref{le-3.4.eq1} follows.
$\qquad \Box$\end{pf}

Now let us start to prove Theorem \ref{matrixstability.tm}.

{\it Proof of Theorem \ref{matrixstability.tm}.}\quad Let $N\ge 1$ be a sufficiently large integer
determined later, $n\in N{\mathbb Z}^d$,  the multiplication operator
$\Psi_n^N$ be as in (\ref{multiplicationoperator.def}), and the
 truncation matrix $A_{N}$ be as in
\eqref{truncationmatrix.def}. By the assumption on the infinite
matrix $A$, there exists a positive constant $C_0$ such that
\begin{equation}\label{matrixstability.tm.pf.eq1}
\|\Psi_n^N c\|_{\ell^p(\Lambda')}
   \le C_0 \|A \Psi_n^N c\|_{\ell^p(\Lambda)}
  \end{equation}
  for any sequence
$c %=(c(\lambda'))_{\lambda'\in\Lambda'}
\in \ell^q(\Lambda')$, $n\in N {\mathbb Z}^d$ and $1\le N\in {\mathbb N}$.
By \eqref{eq-3.1}, \eqref{truncation.eq},  \eqref{eq-3.6},
\eqref{le-3.4.eq1} and \eqref{matrixstability.tm.pf.eq1}, we get
\begin{eqnarray}\label{matrixstability.tm.pf.eq2}
& & \Big(\sum_{n\in N{\mathbb Z}^d} \|\Psi_n^N
c\|_{\ell^p(\Lambda')}^q\Big)^{1/q} \nonumber\\
&\le &  C_0 \Big(\sum_{n\in
N{\mathbb Z}^d} \|A \Psi_n^N
c\|_{\ell^p(\Lambda)}^q\Big)^{1/q}\nonumber\\
& \le & C_0  \Big(\sum_{n\in N{\mathbb Z}^d} \|(A-A_{N}) \Psi_n^N
c\|_{\ell^p(\Lambda)}^q\Big)^{1/q}\nonumber\\
& &
+ C_0 \Big(\sum_{n\in N{\mathbb Z}^d} \| (A_{N} \Psi_n^N-\Psi_n^N
A_N)
c\|_{\ell^p(\Lambda)}^q\Big)^{1/q}\nonumber \\
& &  + C_0 \Big(\sum_{n\in N{\mathbb Z}^d} \| \Psi_n^N (A_N-A)
c\|_{\ell^p(\Lambda)}^q\Big)^{1/q}  + C_0 \Big(\sum_{n\in
N{\mathbb Z}^d} \| \Psi_n^N Ac\|_{\ell^p(\Lambda)}^q\Big)^{1/q}
\nonumber \\
%& \le & C_0 C \|A-A_{N}\|_{\mathcal C} \Big(\sum_{n\in N{\mathbb Z}^d} \|
%\Psi_n^N c\|_{\ell^p(\Lambda)}^q\Big)^{1/q}+ C_0 \Big(\sum_{n\in
%N{\mathbb Z}^d} \|(A_{MN} \Psi_n^N-\Psi_n^N A_{MN}\big)
%c\|_{\ell^p(\Lambda)}^q\Big)^{1/q}\nonumber\\
%& & \quad + C_0 \Big(\sum_{n\in N{\mathbb Z}^d} \|(\Psi_n^N A_{MN}\big)
%c\|_{\ell^p(\Lambda)}^q\Big)^{1/q}\nonumber\\
%& \le & C_0 \big(  \|A-A_{MN}\|_{\mathcal C}+C M^d
%(\|A-A_{\sqrt{N}}\|_{\mathcal C}+N^{-1/2} \|A\|_{\mathcal C}\big)
%\Big(\sum_{n\in N{\mathbb Z}^d} \| \Psi_n^N
%c\|_{\ell^p(\Lambda)}^q\Big)^{1/q}\nonumber\\
%& & + C_0 \Big(\sum_{n\in N{\mathbb Z}^d} \|\Psi_n^N (A-A_{MN})
%c\|_{\ell^p(\Lambda)}^q\Big)^{1/q}+ C_0 \Big(\sum_{n\in N{\mathbb Z}^d}
%\|\Psi_n^N A
%c\|_{\ell^p(\Lambda)}^q\Big)^{1/q}\\
&\le & C_0  C R(\Lambda)^{1/p} R(\Lambda')^{1-1/p} \Big(
\|A-A_{N}\|_{\mathcal C}+ \inf_{0\le s\le
N}\big(\|A-A_{s}\|_{\mathcal C}+\frac{s}{N} \|A\|_{\mathcal
C}\big)\Big) \nonumber\\
& & \qquad \times \Big(\sum_{n\in N{\mathbb Z}^d} \| \Psi_n^N
c\|_{\ell^p(\Lambda')}^q\Big)^{1/q}+ C_0 \Big(\sum_{n\in N{\mathbb Z}^d}
\|\Psi_n^N A c\|_{\ell^p(\Lambda)}^q\Big)^{1/q}\nonumber\\
&\le & C_0  C R(\Lambda)^{1/p} R(\Lambda')^{1-1/p} \inf_{0\le s\le
N}\big(\|A-A_{s}\|_{\mathcal C}+\frac{s}{N} \|A\|_{\mathcal
C}\big) \nonumber\\
& & \qquad \times \Big(\sum_{n\in N{\mathbb Z}^d} \| \Psi_n^N
c\|_{\ell^p(\Lambda')}^q\Big)^{1/q}%\nonumber \\
%& &
+ C_0 \Big(\sum_{n\in N{\mathbb Z}^d} \|\Psi_n^N A
c\|_{\ell^p(\Lambda)}^q\Big)^{1/q}
\end{eqnarray}
where $1\le q<\infty$. Note that for any infinite matrix $A\in
{\mathcal C}(\Lambda, \Lambda')$
\begin{eqnarray*} 0 & \le &
\lim_{N\to \infty}
 \inf_{0\le s\le
N}\big(\|A-A_{s}\|_{\mathcal C}+\frac{s}{N} \|A\|_{\mathcal
C}\big)\\
& \le & \lim_{N\to \infty} \big(\|A-A_{\sqrt{N}}\|_{\mathcal
C}+N^{-1/2} \|A\|_{\mathcal C}\big) =0\end{eqnarray*} by
\eqref{truncation.eq}. Therefore by selecting $N$ sufficiently
large in \eqref{matrixstability.tm.pf.eq2}, we have that
\begin{equation}\label{eq3.11}
\Big(\sum_{n\in N{\mathbb Z}^d} \|\Psi_n^N
c\|_{\ell^p(\Lambda')}^q\Big)^{1/q}\le 2C_0 \Big(\sum_{n\in N{\mathbb Z}^d}
\|\Psi_n^N (Ac)\|_{\ell^p(\Lambda)}^q\Big)^{1/q}.
\end{equation}
By the equivalence of different  norms on a finite-dimensional
space, there exists a positive constant $C$ (that depends on $p,
q, d$ only) such that
\begin{equation}\label{eq3.12a}
C^{-1} (R(\Lambda') N^{d})^{\min(1/p-1/q, 0)} \|\Psi_n^N
c\|_{\ell^q(\Lambda')}\le \|\Psi_n^N c\|_{\ell^p(\Lambda')}%\le C
%(R(\Lambda') N^{d})^{\max(1/p-1/q, 0)} \|\Psi_n^N
%c\|_{\ell^q(\Lambda')}
\end{equation}
and
\begin{equation}\label{eq3.12b}
\|\Psi_n^N Ac\|_{\ell^p(\Lambda)}\le C (R(\Lambda)
N^{d})^{\max(1/p-1/q, 0)} \|\Psi_n^N Ac\|_{\ell^q(\Lambda)}
\end{equation}
 hold for all sequences $c\in \ell^q(\Lambda')$, $n\in N{\mathbb Z}^d$
and $1\le N\in {\mathbb N}$. Therefore combining \eqref{eq3.11},
\eqref{eq3.12a} and \eqref{eq3.12b}, we conclude that
\begin{equation}
\|c\|_{\ell^q(\Lambda')}\le C (R(\Lambda') N^{d})^{-\min (1/p-1/q,
0)} (R(\Lambda) N^{d})^{\max (1/p-1/q, 0)}
\|Ac\|_{\ell^q(\Lambda)}
\end{equation}
for any $c\in \ell^q(\Lambda')$, and the conclusion for $1\le
q<\infty$ follows.

The conclusion for $q=\infty$   can be proved by similar argument. We omit the details here.
$\qquad \Box$

\bigskip

%%==================================================================
\section{Stability for localized synthesis operators}
\label{priesz.section}
%%==================================================================

In this section, we consider the stability of the synthesis
operator
\begin{equation}\label{sphi.def}
S_\Phi:\ \ell^p(\Lambda)\ni (c(\lambda))_{\lambda\in
\Lambda}\longmapsto \sum_{\lambda\in \Lambda}
c(\lambda)\phi_\lambda \in V_p(\Phi, \Lambda)
\end{equation}
associated with a family
 $\Phi= \{ \phi_{\lambda}:
\lambda\in \Lambda \}$ of functions  on ${\mathbb R}^d$, where
\begin{equation} \label{synthesis.tm.eq3}
V_p(\Phi,\Lambda):=\left\{\sum_{\lambda\in \Lambda} c(\lambda)
\phi_\lambda: \ (c_\lambda)_{\lambda\in \Lambda}\in
\ell^p(\Lambda)\right\}, \ 1\le p\le \infty
\end{equation}
(\cite{sunaicm}).
The synthesis operator $S_\Phi$ appears in the study of
 spline approximation and operator  approximation (\cite{devore, schumakerbook}),
 wavelet analysis
(\cite{chuibook, daubechiesbook, mallatbook, meyerbook}), Gabor analysis (\cite{grochenigbook}) and
sampling (\cite{akramgrochenig01, sunsiam}), while one of basic
assumptions for the synthesis operator $S_\Phi$ is  the $\ell^p$-stability,
 i.e., there exists a positive constant $C$ such that
\begin{equation}
C^{-1} \|c\|_{\ell^p(\Lambda)} \le \|S_\Phi c \|_p\le C \| c\|_{\ell^p(\Lambda)} \quad {\rm for\ all}  \
c\in \ell^p(\Lambda).
\end{equation}

 The main results of this section  are Theorem
\ref{synthesis.tm} (a generalization of Theorem
\ref{prieszpreliminary.thm}) about equivalence of the
$\ell^p$-stability of the synthesis operator $S_\Phi$ for
different $1\le p\le\infty$, and  Corollary \ref{synthesisdual.cr}
about  well localization of
  the inverse of  the synthesis operator $S_\Phi$.
% and also a partial answer to a question .

\begin{theorem} \label{synthesis.tm}
Let $\Lambda$ be a relatively-separated subset of ${\mathbb R}^d$,
$\Phi=\{\phi_\lambda, \lambda\in \Lambda\}$ be a family of
functions with the property that
\begin{equation}
\label{synthesis.tm.eq1}\Big\|\sup_{\lambda\in
\Lambda}|\phi_\lambda(\cdot+\lambda)|\ \Big\|_{{\mathcal
W}_1}<\infty
\end{equation}
and
\begin{equation}\label{synthesis.tm.eq2}\lim_{\delta\to 0}\Big\|
\sup_{\lambda\in \Lambda}|\omega_\delta(\phi_\lambda)(\cdot+\lambda)|\
\Big\|_{{\mathcal W}_1}=0.\end{equation}
 If the synthesis operator
$S_\Phi$ in \eqref{sphi.def} has $\ell^p$-stability for some $1\le p\le
\infty$, then it has $\ell^q$-stability
for any $1\le q\le \infty$.
\end{theorem}

For $\Phi=\{\phi_n(\cdot-j)\}_{1\le n\le N, j\in {\mathbb Z}^d}$ generated
by integer shifts of finitely many functions $\phi_1, \ldots,
\phi_N$, we have  the following corollary for the synthesis
operator $S_\Phi$ associated with $\Phi$. Here in the statement of
the following result, we do not include the regularity condition
\eqref{synthesis.tm.eq2} because $\lim_{\delta\to
0}\|\omega_\delta(f)\|_{{\mathcal W}_1}=0$ for any continuous
function $f$ in the Wiener amalgam space ${\mathcal W}_1$
(\cite{akramgrochenig01}).

\begin{corollary}
Let   $\phi_1, \ldots, \phi_N$ be continuous functions in the
Wiener amalgam space ${\mathcal W}_1$, and for $1\le p\le \infty$
define
$$
V_p(\phi_1, \ldots, \phi_N) := \left\{\sum_{n=1}^N
\sum_{j\in{\mathbb Z}^d} c_n(j) \phi_n(\cdot-j): \ (c_n(j))
%_{1\le n\le N,
%j\in {\mathbb Z}^d}
 \in (\ell^p({\mathbb Z}^d))^N\right\}.
$$  If the synthesis operator
$L_{\phi_1, \cdots, \phi_n}:\  (\ell^p({\mathbb Z}^d))^N\longmapsto
V_p(\phi_1, \ldots, \phi_N)$ defined by
$$
L_{\phi_1,\cdots,\phi_n}:\ (c_n(j))_{1\le n\le N, j\in {\mathbb Z}^d}
\longmapsto \sum_{n=1}^N \sum_{j\in{\mathbb Z}^d} c_n(j) \phi_n(\cdot-j)
$$
has  $\ell^p$-stability for some $p\in[1,\infty]$, i.e., there exists a
positive constant $C$ such that
$$ C^{-1} \|c\|_{(\ell^p({\mathbb Z}^d))^N}\le \|L_{\phi_1, \ldots,
\phi_N}c\|_{p}\le C \|c\|_{(\ell^p({\mathbb Z}^d))^N} \quad {\rm for \
all} \ c\in (\ell^p({\mathbb Z}^d))^N, $$  then the synthesis operator
$L_{\phi_1,\cdots,\phi_n}$ has $\ell^q$-stability for any
$q\in[1,\infty]$.
\end{corollary}

The result in the above corollary is established in
\cite{jiamicchelli} under the weak assumption that $\phi_1,
\ldots, \phi_N\in {\mathcal L}^\infty$.

\bigskip

 Note that  the synthesis operator $S_\Phi$  has  $\ell^2$-stability if
and only if the matrix $A=\big (a(\lambda, \lambda')\big)_{\lambda, \lambda'\in \Lambda}$
 has $\ell^2$-stability
 where $a(\lambda, \lambda')=\int_{{\mathbb R}^d} \phi_\lambda(x) \phi_{\lambda'}(x)dx$
for $\lambda, \lambda'\in \Lambda$. This observation  together with the equivalence  in Theorem \ref{synthesis.tm}
 for the synthesis operator $S_\Phi$
 and  the Wiener's lemma in \cite{sjostrand}
 for the Sj\"ostrand class of infinite matrices leads to the following
 result.

\begin{corollary}\label{synthesisdual.cr}
Let $1\le p\le\infty$,  $\Lambda$ be a relatively-separated subset of ${\mathbb R}^d$,
$\Phi=\{\phi_\lambda, \lambda\in \Lambda\}$ satisfy
\eqref{synthesis.tm.eq1} and
\eqref{synthesis.tm.eq2}. If the synthesis operator $S_\Phi$ has  $\ell^p$-stability,
then there exists another family $\tilde \Phi=\{\tilde
\phi_\lambda, \lambda\in \Lambda\}$  functions  satisfying
\eqref{synthesis.tm.eq1} and
\eqref{synthesis.tm.eq2} such that the inverse of the synthesis operator $S_\Phi$ is given by
$$(S_\Phi)^{-1} f= \Big(\int_{{\mathbb R}^d} f(x) \tilde
\phi_\lambda(x)dx\Big)_{\lambda\in \Lambda}\quad {\rm for\ all} \ f\in
V_p(\Phi, \Lambda).$$
\end{corollary}

The conclusion in the above corollary with $p=2$ is established in
 \cite{sunaicm} without the regularity assumption
 \eqref{synthesis.tm.eq2}. The conclusion in the above corollary
 for general $1\le p\le \infty$ gives a partial answer to a problem in \cite[Remark 5.3]{sunaicm}.

\bigskip

To prove Theorem \ref{synthesis.tm}, we recall a result in \cite{sunaicm}.

\begin{lemma}\label{synthesis.lm1}
Let $1\le p\le \infty$, $\Lambda$ be a relatively-separated subset of
${\mathbb R}^d$,  $\Phi=\{\phi_\lambda, \lambda\in \Lambda\}$  satisfy \eqref{synthesis.tm.eq1}.
Then there exists a positive constant $C$
(that depends on $d$ and $p$ only) such that
\begin{equation}
\| S_\Phi c\|_{p}\le  C  R(\Lambda)^{1-1/p} \Big\|\sup_{\lambda\in \Lambda} |\phi_\lambda(\cdot+\lambda)|\ \Big\|_{{\mathcal W}_1} \|c\|_{\ell^p(\Lambda)}.
\end{equation}
\end{lemma}

Now we start to prove  Theorem \ref{synthesis.tm}.

%\begin{pf}
{\it Proof of Theorem \ref{synthesis.tm}.}\quad
Let $1\le p,q\le \infty$. By the $\ell^p$-stablity of the synthesis operator $S_\Phi$, there exists a positive constant
$C_0$ such that
\begin{equation}\label{synthesis.tm.pf.eq1}
\|c\|_{{\ell^p(\Lambda)}}\le C_0 \|S_\Phi c\|_p
\quad {\rm for\ all} \ c\in \ell^p(\Lambda).
\end{equation}

For $1\le n\in {\mathbb N}$, define the operator $P_n$ on $L^p$ by
\begin{equation} \label{synthesis.tm.pf.eq1b}
P_nf(x)=2^{nd} \sum_{\lambda'\in 2^{-n}{\mathbb Z}^d}\phi_0(2^n(x-\lambda'))\times
 \int_{{\mathbb R}^d} f(y)  \phi_0(2^n(y-\lambda')) dy,\ \ f\in L^p({\mathbb R}^d)\end{equation}
where $\phi_0$ be the characteristic function on $[0,1)^d$,
and let $\Phi_n=\{P_n\phi_\lambda,\lambda\in \Lambda\}$.
Then
\begin{equation} \label{synthesis.tm.pf.eq2}
|\phi_\lambda(x)-P_n\phi_\lambda(x)|\le
\omega_{2^{-n}}(\phi_\lambda)(x) \quad {\rm for\ all} \ x\in {\mathbb R}^d \ {\rm and}  \ \lambda\in \Lambda.
\end{equation}
From \eqref{synthesis.tm.pf.eq1}, \eqref{synthesis.tm.pf.eq2} and Lemma \ref{synthesis.lm1} it follows that
\begin{eqnarray} \label{synthesis.tm.pf.eq3}
\|S_\Phi c\|_p & \le &  \| S_{\Phi_n} c\|_p+\|S_{\Phi-\Phi_n}c\|_p\nonumber\\
&\le &
\|S_{\Phi_n}c\|_p +C R(\Lambda)^{1-1/p} \Big\|\sup_{\lambda\in \Lambda}
\omega_{2^{-n}}(\phi_\lambda)(\cdot+\lambda) \Big\|_{{\mathcal W}_1}
\|c\|_{\ell^p(\Lambda)}.
\end{eqnarray}
Combining \eqref{synthesis.tm.eq2} and \eqref{synthesis.tm.pf.eq3} leads to
the existence of a sufficiently large integer $n_0$ such that
\begin{equation} \label{synthesis.tm.pf.eq4}
\|c\|_{{\ell^p(\Lambda)}}\le 2 C_0 \|S_{\Phi_{n_0}} c\|_p
\quad {\rm for\ all} \ c\in \ell^p(\Lambda).
\end{equation}

Define  $A_{n_0}=(a_{n_0}(\lambda', \lambda))_{ \lambda'\in
2^{-n_0}{\mathbb Z}^d, \lambda\in \Lambda}$ by
 \begin{equation} \label{synthesis.tm.pf.eq4b}
 a_{n_0}(\lambda', \lambda)= 2^{n_0d} \int_{{\mathbb R}^d} \phi_\lambda (y) \phi_0(2^{n_0}(y-\lambda')) dy.\end{equation}
Since
$$P_{n_0}\phi_\lambda= \sum_{\lambda'\in 2^{-n_0}{\mathbb Z}^d} a_{n_0}(\lambda', \lambda)  \phi_0(2^{n_0}(\cdot-\lambda')),$$
and
\begin{equation} \label{synthesis.tm.pf.eq4c}
\Big\| \sum_{\lambda'\in 2^{-n_0}{\mathbb Z}^d} a(\lambda') \phi_0(2^{n_0}(\cdot-\lambda'))\Big\|_p=
2^{-n_0d/p} \|a\|_{\ell^p(2^{-n_0}{\mathbb Z}^d)}\end{equation}
for any $a=(a(\lambda'))_{\lambda'\in 2^{-n_0}{\mathbb Z}^d}\in \ell^p(2^{-n_0}{\mathbb Z}^d)$,
 the equation \eqref{synthesis.tm.pf.eq4} can be rewritten
in the following matrix formulation:
 \begin{equation} \label{synthesis.tm.pf.eq5}
\|c\|_{\ell^p(\Lambda)}\le 2C_0 2^{-n_0d/p}
\|A_{n_0}c\|_{\ell^p(2^{-n_0}{\mathbb Z}^d)} \quad {\rm for\ all} \ c\in \ell^p(\Lambda).
 \end{equation}

By \eqref{synthesis.tm.eq1}, it holds that
 \begin{eqnarray*} & & \sum_{j\in {\mathbb Z}^d}
\sup_{\lambda'\in
 2^{-n_0}{\mathbb Z}^d, \lambda\in \Lambda} |a_{n_0}(\lambda', \lambda)| \chi_{j+[0,1)^d}(\lambda'-\lambda)\\
 & \le &
2^{n_0d}
\sum_{j\in {\mathbb Z}^d} \sup_{\lambda'\in
 2^{-n_0}{\mathbb Z}^d, \lambda\in \Lambda} \chi_{j+[0,1)^d}(\lambda'-\lambda)\times
\int_{{\mathbb R}^d} h(y-\lambda) \phi_0(2^{n_0}(y-\lambda')) dy\nonumber\\
& \le & \sum_{j\in {\mathbb Z}^d}
\sup_{y\in j+[0,2)^d} h(y) \le  2^d \|h\|_{{\mathcal W}_1}
<\infty
\end{eqnarray*}
where $h(x)=\sup_{\lambda\in \Lambda} |\phi_\lambda(x+\lambda)|$, which means  that
the infinite matrix $A_{n_0}$ in  \eqref{synthesis.tm.pf.eq4b} belongs to the Sj\"ostrand class ${\mathcal C}(2^{-n_0}{\mathbb Z}^d, \Lambda)$,
\begin{equation} \label{synthesis.tm.pf.eq9}
A_{n_0}\in {\mathcal C}(2^{-n_0}{\mathbb Z}^d, \Lambda).\end{equation}

By \eqref{synthesis.tm.pf.eq5}, \eqref{synthesis.tm.pf.eq9} and
 Theorem \ref{matrixstability.tm},
the infinite matrix $A_{n_0}$ has the  $\ell^q$-stability, i.e.,
there exists a positive constant $C_1$ such that
\begin{equation} \label{synthesis.tm.pf.eq10}
\|c\|_{\ell^q(\Lambda)} \le C_1 \|A_{n_0}c\|_{\ell^q(2^{-n_0}{\mathbb Z}^d)} \quad {\rm for \ all} \ c\in \ell^q(\Lambda).
\end{equation}
For any $c=(c(\lambda))_{\lambda\in \Lambda}\in \ell^q(\Lambda)$,
 \begin{equation} \label{synthesis.tm.pf.eq11}
\|A_{n_0}c\|_{\ell^q(2^{-n_0}{\mathbb Z}^d)}= 2^{n_0d/q}
\Big\| \int_{{\mathbb R}^d} K_{n_0}(\cdot, y) (S_\Phi c)(y)dy\Big \|_q
\le  %2^{n_0d/p}
% (\langle f, 2^{n_0d}
%\phi_0(2^{n_0}-k)\rangle)_{k\in {\mathbb Z}^d}\|_{\ell^p({\mathbb Z}^d)} \le
2^{n_0d/q} \|S_\Phi c\|_{q}
\end{equation}
by
\eqref{synthesis.tm.pf.eq4c}, where
$K_{n_0}(x,y)=2^{n_0d}
\sum_{\lambda'\in 2^{-n_0}{\mathbb Z}^d} \phi_0(2^{n_0}(x-\lambda')) \phi_0(2^{n_0}(y-\lambda'))$.
The $\ell^q$-stability of the synthesis operator
$S_\Phi$ then follows from  \eqref{synthesis.tm.pf.eq10} and
\eqref{synthesis.tm.pf.eq11}.

\bigskip

\section{$L^p$-stability for localized integral operators}
\label{localizedintegraloperator.section}

In this section, we consider the $L^p$-stability of  integral
operators
\begin{equation}\label{lio.def}
Tf(x):= \int_{{\mathbb R}^d} K_T(x, y) f(y)dy, \quad   f\in L^p({\mathbb R}^d)
\end{equation}
whose kernels $K_T$  are enveloped by  convolution kernels with certain decay at infinity, i.e.,
\begin{equation}
|K_T(x,y)|\le h(x-y)\ {\rm for \ all} \ x, y\in {\mathbb R}^d
\end{equation}
 where $h$ is a  function  in the Wiener amalgam space ${\mathcal W}_1$ (\cite{barnes,
 brandenburg, jorgens,
 Kurbatov01, Kurbatov, sunacha}).
 Examples of the integral operators of the form  \eqref{lio.def}
 include  projection operators on  wavelet spaces (\cite{chuisun, cohen96,
coifmanappear, daubechiesbook, jiamicchelli, sunaicm}),  frame operators associated with  Gabor systems in the time-frequency space (\cite{balan, balanchl04, grochenigbook, grochenigl03}), and
reconstruction  operators in sampling theory  (\cite{akramgrochenig01,  sunsiam, sunaicm}).

An integral operator $T$ with kernel $K_T$ enveloped by a
convolution kernel in the Wiener amalgam space defines  a bounded
operator on $L^p$.  The above class ${\mathcal C}_1$ of  localized
integral operators  is a non-unital algebra. The new  algebra
 $${\mathcal I}{\mathcal C}_1=\{ \lambda I+ T: \ \lambda\in {\mathbb C}, T\in {\mathcal C}_1\}$$
 obtained by adding the identity operator $I$ on $L^p$
 to that algebra ${\mathcal C}_1$ is a  unital Banach subalgebra of ${\mathcal B}(L^p), 1\le p\le \infty$ (\cite{sunacha}).

 In this section, we discuss the $L^p$-stability of the  localized integral
 operators in
${\mathcal I}{\mathcal C}_1$ with additional regularity on
kernels. The main results of this section are Theorem
\ref{operatorstability.tm} %concerning the equivalence of
%$L^p$-stability of localized integral operator
(a slight
generalization of Theorem \ref{liopreliminary.thm}), and Corollary
\ref{liowiener.cr} concerning the well localization of the inverse
of a localized integral operator.

%The  Wiener's lemma for integral (pseudo-differential) operators
%is also an important chapter of the theory of inverse-closed
%subalgebras (\cite{barnes, baskakov, boulkhemair1, boulkhemair2,
%  farellstrohmer,
%grochenigl03, grocheniglappear, domar56,  herau, Kurbatov01,
%sjostrand, sjostrand2, toft1, toft2, toft3} and also \cite{
%grochenigappear, Kurbatov} for a historial review).
%
%
%
%
%In \cite{sun-wie}, the Wiener's lemma for localized integral operator is proved.
%The kernel function $K_T$ belongs to the class $W_{p,u}^\alpha$ and the integral operator
%$T$ on $L^2({\mathbb R}^d)$ is defined by
%$$
%(Tf)(x)=\int K_T(x,y) f(y) dy.
%$$
%It is showed in \cite{sun-wie} that if $I+T$ has an inverse in
%${\mathcal B}^2$(the space of bounded operators in $L^2({\mathbb R}^d)),$
%then the inverse  belongs to ${\sl IW}_{p,u}^\alpha.$
%
%
%
%
%In this section, we consider the equivalence  of $L^p$-stability of localized integral operators.

\begin{theorem}\label{operatorstability.tm}
 Let $0<\alpha\le 1$,  $D$ be a positive constant, and  $T$ be an integral operator of the form \eqref{lio.def} with
 its kernel $K_T$  satisfying
\begin{equation}\label{operatorstability.tm.eq1}
\Big\|\sup_{y\in {\mathbb R}^d}|K_T(y, \cdot+y)| \Big\|_{ {\mathcal W}_1}\le D\end{equation}
 and
 \begin{equation} \label{operatorstability.tm.eq2}
\Big\|\sup_{y\in {\mathbb R}^d} \omega_\delta(K_T)(y, \cdot+y)\Big\|_{{\mathcal W}_1}\le D \delta^\alpha \ {\rm for \ all} \ \delta\in (0,1). \end{equation}
If $I+T$ has $L^p$-stability for some $1\le p\le \infty$,
then it has $L^q$-stability for all $1\le q\le \infty$.
\end{theorem}

The above result  for $p=2$ follows from the Wiener's lemma for localized integral operators (\cite{sunacha}).
Applying the $L^p$ equivalence in Theorem \ref{operatorstability.tm} for different $p$, we can extend
  the Wiener's lemma in \cite{sunacha} to $p\ne 2$.

  \begin{corollary}\label{liowiener.cr}
   Let  $1\le p\le \infty, 0\ne \lambda\in {\mathbb C}$,
   and  $T$ be an integral operator with its kernel $K_T$  satisfying
\eqref{operatorstability.tm.eq1} and \eqref{operatorstability.tm.eq2}.
If $\lambda I+T$ has bounded inverse
on  $L^p$, then $(\lambda I+T)^{-1}= \lambda^{-1} I+ \tilde T$
for some integral operator $\tilde T$ with kernel satisfying
\eqref{operatorstability.tm.eq1} and \eqref{operatorstability.tm.eq2}.
  \end{corollary}

  Recall that any integral operator having its kernel satisfying
  \eqref{operatorstability.tm.eq1} and
  \eqref{operatorstability.tm.eq2} does not have bounded inverse
  in $L^p, 1\le p<\infty$ (\cite{sunacha}). Then from Corollary
  \ref{liowiener.cr} we have the
  following result to spectra of localized integral operators on $L^p$.

\begin{corollary}
   Let   $T$ be an integral operator with its kernel $K_T$  satisfying
\eqref{operatorstability.tm.eq1} and
\eqref{operatorstability.tm.eq2}. Then
 \begin{equation}
 \sigma_p(T)=\sigma_q(T)\quad {\rm for \ all} \ 1\le p, q<\infty,
 \end{equation}
 where $\sigma_p(T)$ denotes the spectrum of the operator $T$ on $L^p$.
\end{corollary}

 Now we start to prove Theorem \ref{operatorstability.tm}.

%\begin{pf}
{\it Proof of Theorem \ref{operatorstability.tm}.}\quad
By the  $L^p$-stability of the operator $I+T$, there exists a positive constant $C_0$ such that
\begin{equation}\label{operatorstability.tm.pf.eq1}
\|f\|_p\le C_0\|(I+T)f\|_p \quad {\rm for  \ all} \ f\in L^p.
\end{equation}

For $1\le n\in {\mathbb N}$, let $T_n=P_nTP_n$ with kernel $K_{T_n}$
where $P_n$ is given in  \eqref{synthesis.tm.pf.eq1b}.
Then
\begin{eqnarray} \label{operatorstability.tm.pf.eq1b}
K_{T_n}(x,y) %& = & \int_{{\mathbb R}^d} \int_{{\mathbb R}^d} K_n(x, s) K_T(s, t) K_n(t, y) ds dt\\
& = & \sum_{\lambda, \lambda'\in 2^{-n} {\mathbb Z}^d} a_n(\lambda, \lambda')   \phi_0(2^n(x-\lambda))\phi_0(2^n(y-\lambda'))
\end{eqnarray}
and
\begin{equation} \label{operatorstability.tm.pf.eq3}
|K_{T_n}(x, y)-K_T(x,y)|\le C \omega_{2^{-n}}(K_T)(x,y) \quad {\rm
for \ all} \ x, y\in {\mathbb R}^d,\end{equation} where $\phi_0$ is the
characteristic function on $[0,1)^d$ and
\begin{equation} \label{operatorstability.tm.pf.eq3b}
a_n(\lambda, \lambda')= 2^{2nd} \int_{{\mathbb R}^d} \int_{{\mathbb R}^d} \phi_0(2^n(s-\lambda))
 K_T(s, t) \phi_0(2^n(t-\lambda'))ds dt\end{equation}
for $\lambda, \lambda'\in 2^{-n}{\mathbb Z}^d$. Therefore we have from
\eqref{operatorstability.tm.eq2} and
\eqref{operatorstability.tm.pf.eq3} that for any $ f\in L^r$ with
$1\le r\le \infty$,
\begin{eqnarray} \label{operatorstability.tm.pf.eq4}
 \|(T-T_n)f\|_r
 &\le &
C \Big \|\sup_{y\in {\mathbb R}^d}  \omega_{2^{-n}}(K_T) (y,
\cdot+y)\Big\|_{{\mathcal W}_1}\|f\|_r\nonumber\\
 &\le & C 2^{-n\alpha}
\|f\|_r.
\end{eqnarray}

By  \eqref{operatorstability.tm.eq2},
 \eqref{operatorstability.tm.pf.eq1}  and \eqref{operatorstability.tm.pf.eq4},
there exists a sufficiently large  integer $n_0$
such that for all $n\ge n_0$,
\begin{equation}\label{operatorstability.tm.pf.eq5}
\|f\|_p\le 2C_0\|(I+T_n)f\|_p \quad {\rm for  \ all} \ f\in L^p.
\end{equation}

Let
\begin{equation}
A_n:=(a_n(\lambda, \lambda'))_{\lambda, \lambda'\in 2^{-n} {\mathbb Z}^d}
\end{equation}
where  $a_n(\lambda, \lambda'), \lambda, \lambda'\in 2^{-n}{\mathbb Z}^d$,
are given in \eqref{operatorstability.tm.pf.eq3b}.
Applying  \eqref{operatorstability.tm.pf.eq5} to
 $$f_n:=\sum_{\lambda\in  2^{-n} {\mathbb Z}^d} c_n(\lambda)
\phi_0(2^n(\cdot-\lambda))\quad {\rm  with} \ (c_n(\lambda))_{\lambda\in  2^{-n}{\mathbb Z}^d}\in \ell^p(2^{-n}{\mathbb Z}^d),$$
and noting
\begin{equation} \label{operatorstability.tm.pf.eq5b}
\| f_n\|_p= 2^{-nd/p} \|c_n\|_{\ell^p(2^{-n}{\mathbb Z}^d)}
\end{equation}
and
\begin{equation} \label{operatorstability.tm.pf.eq6}
\| (I+T_n)f_n\|_p= 2^{-nd/p} \| (I+ 2^{-nd} A_n)c_n\|_{\ell^p(2^{-n}{\mathbb Z}^d)},
\end{equation}
we obtain the uniform  $\ell^p$-stability of the matrix $I+2^{-nd}
A_n$, i.e.,
\begin{equation} \label{operatorstability.tm.pf.eq8}
\|c_n\|_{\ell^p(2^{-n}{\mathbb Z}^d)}\le 2C_0 \| (I+ 2^{-nd}
A_n)c_n\|_{\ell^p(2^{-n}{\mathbb Z}^d)} \end{equation} holds for  any $
c_n\in \ell^p(2^{-n}{\mathbb Z}^d)$ and $n\ge n_0$.

Define
\begin{equation} \label{operatorstability.tm.pf.eq9}
A_{n,s}= (a_{n, s}(\lambda, \lambda'))_{\lambda, \lambda'\in 2^{-n} {\mathbb Z}^d}\end{equation}
where
$$a_{n,s}(\lambda, \lambda')=\left\{\begin{array}{ll}
a_{n}(\lambda, \lambda') & {\rm  if} \ \|\lambda-\lambda'\|_\infty<s,\\
0 & {\rm otherwise}.\end{array}\right.$$ Then for $s\ge 0$,
\begin{eqnarray*} %\label{operatorstability.tm.pf.eq10}
\| A_n-A_{n,s}\|_{\mathcal C} & \le &
\sum_{j\in {\mathbb Z}^d \ {\rm with} \ \|j\|_\infty\ge s-1}
\sup_{\lambda, \lambda'\in 2^{-n}{\mathbb Z}^d} |a_n(\lambda, \lambda')| \chi_{j+[0,1)^d}(\lambda-\lambda')\nonumber\\
& \le & 2^{2nd}
\sum_{j\in {\mathbb Z}^d \ {\rm with} \ \|j\|_\infty\ge s-1} \sup_{\lambda, \lambda'\in 2^{-n}{\mathbb Z}^d}  \chi_{j+[0,1)^d}(\lambda-\lambda')
\nonumber\\
& & \qquad\times
\int_{2^{-n}[0,1)^d} \int_{2^{-n}[0,1)^d} |K(\lambda+s, \lambda'+t)| dsdt \nonumber\\
& \le &
\sum_{j\in {\mathbb Z}^d \ {\rm with} \ \|j\|_\infty\ge s-1} \sup_{x\in j+[-1, 2)^d} |h(x)|\nonumber\\
& \le & 3^d  \sum_{j\in {\mathbb Z}^d \ {\rm with} \ \|j\|_\infty\ge s-3}
\sup_{x\in j+[0, 1)^d} |h(x)|
\end{eqnarray*}
where  $h(x)=\sup_{y\in {\mathbb R}^d} |K_T(x+y, y)|$. Thus
\begin{eqnarray} \label{operatorstability.tm.pf.eq10}
& & \inf_{0\le s\le N}\big(\|A_n-A_{n,s}\|_{\mathcal
C}+\frac{s}{N} \|A_n\|_{\mathcal C}\big)\nonumber\\
& \le & \|A_n-A_{n,\sqrt{N}}\|_{\mathcal
C}+N^{-1/2} \|A_n\|_{\mathcal C}\nonumber\\
%&\le  &  \sum_{j\in {\mathbb Z}^d \ {\rm with} \ \|j\|_\infty\ge
%\sqrt{N}-1} \sup_{x\in j+[-1, 2)^d} |h(x)|\nonumber\\
%& & + N^{-1/2} \sum_{j\in {\mathbb Z}^d}
%\sup_{x\in j+[-1, 2)^d} |h(x)|\nonumber\\
 &
\le & C\Big(\sum_{j\in {\mathbb Z}^d \ {\rm with} \ \|j\|_\infty\ge
\sqrt{N}-3} \sup_{x\in j+[0,1)^d} |h(x)|\nonumber\\
& & + N^{-1/2} \sum_{j\in {\mathbb Z}^d} \sup_{x\in j+[0, 1)^d}
|h(x)|\Big)
\end{eqnarray}

Let $N$ be a sufficiently large integer chosen later and
the multiplication operator
$\Psi_j^N$ be as in the proof of Theorem \ref{matrixstability.tm}. Then   for $1\le q\le \infty$,
using
 the similar  argument  in the proof of  Theorem \ref{matrixstability.tm}, we obtain
 from \eqref{truncation.eq}, \eqref{operatorstability.tm.pf.eq8},  \eqref{operatorstability.tm.pf.eq10}, and Lemmas \ref{lpboundedness.lm}, \ref{multiplicationoperation.lm} and \ref{le-3.4}
that
 \begin{eqnarray}\label{operatorstability.tm.pf.eq11}
& & \left\| \left( \|\Psi_j^N
c_n\|_{\ell^p( 2^{-n} {\mathbb Z}^d)}\right)_{j\in N{\mathbb Z}^d}\right\|_{\ell^q(N{\mathbb Z}^d)} \nonumber\\
&\le &  2C_0 \left\| \left(  \| (I+ 2^{-nd}A_n) \Psi_j^N
c_n\|_{\ell^p(2^{-n}{\mathbb Z}^d)}\right)_{j\in N{\mathbb Z}^d}\right\|_{\ell^q(N{\mathbb Z}^d)} \nonumber\\
& \le & 2^{-nd+1}C_0   \left\| \left(  \|(A_n-A_{n,N}) \Psi_j^N
c_n\|_{\ell^p(2^{-n}{\mathbb Z}^d)}\right)_{j\in N{\mathbb Z}^d}\right\|_{\ell^q(N{\mathbb Z}^d)} \nonumber\\
& &  + 2^{-nd+1}C_0 \left\|\left( \| (A_{n, N} \Psi_j^N-\Psi_j^N
A_{n,N})
c_n\|_{\ell^p(2^{-n}{\mathbb Z}^d)}\right)_{j\in N{\mathbb Z}^d}\right\|_{\ell^q(N{\mathbb Z}^d)} \nonumber \\
& &  + 2^{-nd+1}C_0 \left\| \left( \| \Psi_j^N (A_{n,N}-A_n)
c_n\|_{\ell^p(2^{-n}{\mathbb Z}^d)}\right)_{j\in N{\mathbb Z}^d}\right\|_{\ell^q(N{\mathbb Z}^d)} \nonumber\\
& &  + 2C_0 \left\| \left( \| \Psi_j^N (I+ 2^{-nd}
A_n)c_n)\|_{\ell^p(\Lambda)}\right)_{j\in
N{\mathbb Z}^d}\right\|_{\ell^q(N{\mathbb Z}^d)}
\nonumber \\
%& \le & C_0 C \|A-A_{N}\|_{\mathcal C} \Big(\sum_{n\in N{\mathbb Z}^d} \|
%\Psi_n^N c\|_{\ell^p(\Lambda)}^q\Big)^{1/q}+ C_0 \Big(\sum_{n\in
%N{\mathbb Z}^d} \|(A_{MN} \Psi_n^N-\Psi_n^N A_{MN}\big)
%c\|_{\ell^p(\Lambda)}^q\Big)^{1/q}\nonumber\\
%& & \quad + C_0 \Big(\sum_{n\in N{\mathbb Z}^d} \|(\Psi_n^N A_{MN}\big)
%c\|_{\ell^p(\Lambda)}^q\Big)^{1/q}\nonumber\\
%& \le & C_0 \big(  \|A-A_{MN}\|_{\mathcal C}+C M^d
%(\|A-A_{\sqrt{N}}\|_{\mathcal C}+N^{-1/2} \|A\|_{\mathcal C}\big)
%\Big(\sum_{n\in N{\mathbb Z}^d} \| \Psi_n^N
%c\|_{\ell^p(\Lambda)}^q\Big)^{1/q}\nonumber\\
%& & + C_0 \Big(\sum_{n\in N{\mathbb Z}^d} \|\Psi_n^N (A-A_{MN})
%c\|_{\ell^p(\Lambda)}^q\Big)^{1/q}+ C_0 \Big(\sum_{n\in N{\mathbb Z}^d}
%\|\Psi_n^N A
%c\|_{\ell^p(\Lambda)}^q\Big)^{1/q}\\
&\le & C_0  C  \Big( \|A_n-A_{n,N}\|_{\mathcal C}+ \inf_{0\le s\le
N}\big(\|A_n-A_{n,s}\|_{\mathcal C}+\frac{s}{N} \|A_n\|_{\mathcal
C}\big)\Big) \nonumber\\
& & \qquad \times \Big(\sum_{j\in N{\mathbb Z}^d} \| \Psi_j^N
c_n\|_{\ell^p(\Lambda)}^q\Big)^{1/q}+ C_0 \Big(\sum_{j\in N{\mathbb Z}^d}
\|\Psi_j^N A_n c_n\|_{\ell^p(\Lambda)}^q\Big)^{1/q}\nonumber\\
%&\le & CC_0  \inf_{0\le s\le N}\big(\|A_n-A_{n,s}\|_{\mathcal
%C}+\frac{s}{N} \|A_n\|_{\mathcal C}\big) \left\| \left( \|
%\Psi_j^N
%c_n\|_{\ell^p(2^{-nd}{\mathbb Z}^d)}\right)_{j\in N{\mathbb Z}^d}\right\|_{\ell^q(N{\mathbb Z}^d)} \nonumber \\
%& & +2C_0 \left\| \left( \|\Psi_j^N (I+ 2^{-nd}A_n)
%c_n\|_{\ell^p( 2^{-nd} {\mathbb Z}^d)}\right)_{j\in N{\mathbb Z}^d}\right\|_{\ell^q(N{\mathbb Z}^d)}\nonumber\\
%&\le & CC_0  \Big(\sum_{j\in {\mathbb Z}^d \ {\rm with} \ \|j\|_\infty\ge \sqrt{N}-1} \sup_{x\in j+[-1, 2)^d} |h(x)|+
%N^{-1/2} \sum_{j\in {\mathbb Z}^d} \sup_{x\in j+[-1, 2)^d} |h(x)|\Big)\nonumber\\
%& & \qquad \times  \left\| \left(\| \Psi_j^N
%c_n\|_{\ell^p(2^{-nd}{\mathbb Z}^d)}\right)_{j\in N{\mathbb Z}^d}\right\|_{\ell^q(N{\mathbb Z}^d)} \nonumber\\
%& &  +{2 C_0} \left\| \left( \|\Psi_j^N (I+ 2^{-nd}A_n)
%c_n\|_{\ell^p( 2^{-nd} {\mathbb Z}^d)}\right)_{j\in N{\mathbb Z}^d}\right\|_{\ell^q(N{\mathbb Z}^d)}\nonumber\\
&\le & CC_0  \Big(\sum_{j\in {\mathbb Z}^d \ {\rm with} \ \|j\|_\infty\ge
\sqrt{N}-3} \sup_{x\in j+[0,1)^d} |h(x)|\nonumber\\
& & +
N^{-1/2} \sum_{j\in {\mathbb Z}^d} \sup_{x\in j+[0, 1)^d} |h(x)|\Big)%\nonumber\\
%& & \qquad \times
 \left\| \left(\| \Psi_j^N
c_n\|_{\ell^p(2^{-n}{\mathbb Z}^d)}\right)_{j\in N{\mathbb Z}^d}\right\|_{\ell^q(N{\mathbb Z}^d)} \nonumber\\
& &  +{2 C_0} \left\| \left( \|\Psi_j^N (I+ 2^{-nd}A_n)
c_n\|_{\ell^p( 2^{-n} {\mathbb Z}^d)}\right)_{j\in N{\mathbb
Z}^d}\right\|_{\ell^q(N{\mathbb Z}^d)},
 \end{eqnarray}
 where  $h(x)=\sup_{y\in {\mathbb R}^d} |K_T(y, x+y)|$ and
 the uppercase letter $C$ denotes an  absolute constant independent of $n\ge n_0$  and $N\ge 1$
 but may be different at  different occurrences.
By \eqref{operatorstability.tm.eq1} and \eqref{operatorstability.tm.pf.eq11}
 there exists a sufficiently large integer $N_0$ (independent of $n\ge n_0$) such that
\begin{eqnarray} \label{operatorstability.tm.pf.eq12}
 & & \Big\| \Big( \|\Psi_j^{N_0}
c_n\|_{\ell^p( 2^{-n} {\mathbb Z}^d)}\Big)_{j\in
N_0{\mathbb Z}^d}\Big\|_{\ell^q(N_0{\mathbb Z}^d)} \nonumber\\
&\le& 4 C_0 \Big\| \Big(\|\Psi_j^{N_0} (I+ 2^{-nd}A_n)
c_n\|_{\ell^p( 2^{-n} {\mathbb Z}^d)}\Big)_{j\in N_0{\mathbb
Z}^d}\Big\|_{\ell^q(N_0{\mathbb Z}^d)}
\end{eqnarray}
holds for any $c_n\in \ell^q(2^{-n} {\mathbb Z}^d)$ and $n\ge n_0$.

Combining \eqref{eq3.12a}, \eqref{eq3.12b} and
\eqref{operatorstability.tm.pf.eq12} yields
\begin{equation} \label{operatorstability.tm.pf.eq13}
 \|c_n\|_{\ell^q( 2^{-n} {\mathbb Z}^d)}
\le C_1 2^{nd|1/p-1/q|} \| (I+ 2^{-nd}A_n) c_n\|_{\ell^q( 2^{-n}
{\mathbb Z}^d)} \ {\rm for \ all} \ c_n\in \ell^q( 2^{-n} {\mathbb
Z}^d),
\end{equation}
where $C_1$ is a positive constant independent of $n\ge n_0$.
This together with  \eqref{operatorstability.tm.pf.eq5b} and
\eqref{operatorstability.tm.pf.eq6} proves  that
\begin{equation} \label{operatorstability.tm.pf.eq14}
\|f_n\|_q\le  C_1 2^{nd|1/p-1/q|}\|(I+T_n)f_n\|_q
\end{equation}
holds
for any $f_n=\sum_{\lambda\in 2^{-n} {\mathbb Z}^d} c_n(\lambda)\phi_0(2^n(\cdot-\lambda))$ with
$(c_n(\lambda))_{\lambda\in 2^{-n}{\mathbb Z}^d}\in \ell^q(2^{-n}{\mathbb Z}^d)$.

Note that for any $f\in L^q$,
\begin{equation} \label{operatorstability.tm.pf.eq15}
P_nf=  2^{nd} \sum_{\lambda\in 2^{-n}{\mathbb Z}^d}
 \phi_0(2^n(\cdot-\lambda))\times
\int_{{\mathbb R}^d} f(y) \phi_0(2^n(y-\lambda))dy \end{equation} with
\begin{equation} \label{operatorstability.tm.pf.eq16}\Big\|  \Big(2^{nd} \int_{{\mathbb R}^d} f(y) \phi_0(2^n(y-\lambda))dy\Big)_{\lambda\in 2^{-n}{\mathbb Z}^d}\Big\|_{\ell^q(2^{-n}{\mathbb Z}^d)}\le 2^{nd/q} \|f\|_q<\infty.
\end{equation}
Therefore it follows from \eqref{operatorstability.tm.pf.eq14} that
\begin{equation} \label{operatorstability.tm.pf.eq17}
\|P_nf\|_q\le  C_1 2^{nd|1/p-1/q|}\|(I+T_n)P_nf \|_q \quad  {\rm for\ all} \ f\in L^q.
\end{equation}

By \eqref{operatorstability.tm.pf.eq5b}, \eqref{operatorstability.tm.pf.eq15} and
\eqref{operatorstability.tm.pf.eq16},  we have
\begin{equation}\label{operatorstability.tm.pf.eq3-a}
\|P_nf\|_q\le \|f\|_q \ {\rm for \ all} \ f\in L^q.
\end{equation}
This implies that
\begin{equation}\label{operatorstability.tm.pf.eq2}
\|f\|_q\le \|f-P_nf\|_q+\|P_nf\|_q\le 3 \|f\|_q \ {\rm for \ all} \ f\in L^q.
 \end{equation}
Noting that $ P_n^2=P_n$,
$(I-P_n)(I+T)f=(I-P_n)f+(I-P_n)T(I-P_n)f+(I-P_n)T P_n f$
  and
$ P_n(I+T)f=P_n (I+T) P_n f+P_n T (I-P_n)^2 f$,
 we obtain from the second inequality of (\ref{operatorstability.tm.pf.eq2}) that for any $f\in L^q$,
\begin{eqnarray} \label{operatorstability.tm.pf.eq18}
\| (I+T)f\|_q & \ge &  \frac{1}{3} \|(I-P_n) (I+T) f\|_q+\frac{1}{3} \| P_n(I+T)f\|_q\nonumber\\
& \ge & \frac{1}{3} \|(I-P_n)f\|_q+\frac{1}{3} \| P_n(I+T)P_nf\|_q \nonumber\\
& &  -\frac{1}{3} \|(I-P_n)T(I-P_n) f\|_q -\frac{1}{3}
\|(I-P_n)TP_n f\|_q\nonumber\\
& & - \frac{1}{3} \|P_n T (I-P_n)^2 f\|_q.
\end{eqnarray}
We note that $(I-P_n)T$  and $T(I-P_n)$ are integral operators with  their kernel bounded by $\omega_{2^{-n}}(K_T)$
where $K_T$ is the kernel of the integral operator $T$.
Therefore similar to the argument in \eqref{operatorstability.tm.pf.eq4} we have
\begin{equation} \label{operatorstability.tm.pf.eq19}
\| (I-P_n) Tf\|_r+ \|T(I-P_n)f\|_r\le C 2^{-n\alpha} \|f\|_r \ {\rm for \ all } \ f\in L^r
\end{equation}
where $1\le r\le \infty$ and $C$ is a positive constant
independent of $n\ge n_0$.

For those $1\le q\le \infty$ satisfying
$|1/q-1/p|<\alpha/d$, we get from \eqref{operatorstability.tm.pf.eq17},  \eqref{operatorstability.tm.pf.eq3-a},
\eqref{operatorstability.tm.pf.eq2},
\eqref{operatorstability.tm.pf.eq18},
and \eqref{operatorstability.tm.pf.eq19} that
\begin{eqnarray}
\|(I+T)f\|_q &\ge &\Big (\frac{1}{3}-C 2^{-n\alpha}\Big)
\|f-P_nf\|_q \nonumber\\
 &  & + \Big((3C_1)^{-1} 2^{-nd|1/p-1/q|}-C
2^{-n\alpha}\Big)
 \| P_n f\|_q\nonumber\\
 & \ge & \frac{1}{4} \|f-P_nf\|_q+ (4C_1)^{-1} 2^{-nd|1/p-1/q|}\|P_nf\|_q\nonumber\\
 &\ge &
 (4C_1)^{-1} 2^{-nd|1/p-1/q|} \|f\|_q \ {\rm for\ all} \ f\in L^q
\end{eqnarray}
if we let the integer $n$ be chosen to be sufficiently large. This proves
that $I+T$ has  $L^q$-stability if $|1/p-1/q|<\alpha/d$.

Using the above argument iteratively leads to the conclusion that
$I+T$ has $L^q$-stability for all $1\le q\le \infty$.
$\qquad \Box$
%  \end{pf}

The authors thank Professors Akram Aldroubi, Radu Balan,  Ilya
Krishtal,  Deguang Han,  and Romain Tessera for their discussion
and suggestions in preparing the manuscript.

\end{document}